\documentclass[twoside]{gSTA2e}

\usepackage{epstopdf}
\usepackage{subfigure}
\usepackage[bitstream-charter]{mathdesign}
\let\circledS\undefined

\usepackage{dsfont}
\usepackage{amsmath}
\usepackage{amssymb}
\usepackage{dsfont}
\usepackage{mathrsfs}
\usepackage{rotating}
\usepackage{amsthm}
\usepackage{bm}
\usepackage{extarrows}
\usepackage{enumerate}
\usepackage{enumitem}
\usepackage{mathtools}
\usepackage{etoolbox}
\patchcmd{\footnotemark}{\stepcounter{footnote}}{\refstepcounter{footnote}}{}{}
\RequirePackage[OT1]{fontenc}
\RequirePackage{amsthm,amsmath,natbib}
\usepackage{url}
\usepackage{booktabs,multirow,caption}
\usepackage[flushleft]{threeparttable}
\usepackage{rotating}
\usepackage{algorithm}
\usepackage[noend]{algpseudocode}
\usepackage[colorlinks,citecolor=blue,linkcolor=blue,urlcolor=blue]{hyperref}
\usepackage{fancyhdr}
\makeatletter
\usepackage[bottom]{footmisc}
\newcommand{\mpfootnotes}[1][1]{
  \renewcommand{\thempfootnote}{\thefootnote}
  \addtocounter{footnote}{-#1}
  \renewcommand{\footnote}{\stepcounter{footnote}\footnotetext[\value{footnote}]}}
\makeatother

\numberwithin{equation}{section}
\renewcommand{\Pr}{\operatorname{\mathds{P}}}
\newcommand{\E}{\operatorname{\mathds{E}}}
\newcommand{\RR}{\mathds{R}}  \newcommand{\R}{\mathds{R}}
\newcommand{\NN}{\mathds{N}}  \newcommand{\N}{\mathds{N}}
\newcommand{\ZZ}{\mathds{Z}}  
\newcommand{\CC}{\mathds{C}}  
\newcommand{\Var}{\operatorname{\mathsf{Var}}}
\newcommand{\Cov}{\operatorname{\mathsf{Cov}}}
\newcommand{\lead}{\operatorname{lead}}
\newcommand{\Span}{\operatorname{span}}
\newcommand{\sign}{\operatorname{sign}}
\newcommand{\CO}{\mathrm{CO}(\mu;\delta,\beta,\gamma)}
\newcommand{\COq}{\mathrm{CO}(\mu;q)}
\newcommand{\CCC}{\mathcal{C}}
\newcommand{\CH}{\mathcal{H}}
\newcommand{\CX}{\mathcal{X}}
\newcommand{\CS}{\mathcal{S}}
\newcommand{\XA}[2]{\begin{minipage}[A]{40ex}#1\\\vspace{-2ex}\begin{flushright}#2\end{flushright}\end{minipage}}
\newcommand{\XB}[1]{\begin{minipage}[B]{16ex}\begin{flushright}#1\end{flushright}\end{minipage}}
\newcommand{\XC}[1]{\begin{minipage}[C]{30ex}\begin{flushright}#1\end{flushright}\end{minipage}}
\newcommand{\ds}{\displaystyle}

\newcommand{\invsubset}{\mathrel{\reflectbox{\rotatebox[origin=c]{90}{$\subseteq$}}}}
\let\originalleft\left
\let\originalright\right
\renewcommand{\left}{\mathopen{}\mathclose\bgroup\originalleft}
\renewcommand{\right}{\aftergroup\egroup\originalright}
\usepackage{cleveref}
\theoremstyle{plain}
\newtheorem{theorem}{Theorem}[section]
\newtheorem{corollary}[theorem]{Corollary}
\newtheorem{lemma}[theorem]{Lemma}
\newtheorem{proposition}[theorem]{Proposition}

\theoremstyle{definition}
\newtheorem{definition}[theorem]{Definition}
\newtheorem{example}[theorem]{Example}
\newtheorem{application}[theorem]{Application}

\theoremstyle{remark}
\newtheorem{remark}[theorem]{Remark}
\newtheorem{notation}[theorem]{Notation}
\crefname{section}{section}{sections}
\crefname{subsection}{subsection}{subsections}
\crefname{subsubsection}{subsubsection}{subsubsections}
\crefname{paragraph}{paragraph}{paragraphs}
\crefname{theorem}{theorem}{theorems}
\crefname{corollary}{corollary}{corollaries}
\crefname{lemma}{lemma}{lemmas}
\crefname{proposition}{proposition}{propositions}
\crefname{definition}{definition}{definitions}
\crefname{example}{example}{examples}
\crefname{application}{application}{applications}
\crefname{remark}{remark}{remarks}
\crefname{notation}{notation}{notations}
\crefname{algorithm}{algorithm}{algorithms}
\crefname{table}{table}{tables}
\crefname{equation}{equation}{equations}

\begin{document}



\title{Orthogonal polynomials in the Cumulative Ord family and its application to variance bounds}

\author{
\name{Georgios Afendras\textsuperscript{a}$^{\ast}$\thanks{$^\ast$Corresponding author. Email: gafendra@buffalo.edu},
      Narayanaswamy Balakrishnan\textsuperscript{b}$^{\dag}$\thanks{$^\dag$Email: bala@mcmaster.ca}
      and
      Nickos Papadatos\textsuperscript{c}$^{\ddag}$\thanks{$^\ddag$Email: npapadat@math.uoa.gr}}
\affil{\textsuperscript{a}Department of Biostatistics and Jacobs School of Medicine and Biomedical Sciences, The State University of New York at Buffalo,
 Buffalo, NY 14214, USA}
\affil{\textsuperscript{b}Department of Mathematics and Statistics, McMaster University, Hamilton, Ontario, L8S 4K1, Canada}
\affil{\textsuperscript{c}Department of Mathematics, Section of Statistics and O.R.,
 University of Athens, Panepistemiopolis, 157 84 Athens, Greece}
\received{\textcolor{white}{---}}
}

\maketitle

\begin{abstract}
 This article presents and reviews several basic properties of the
 Cumulative Ord family of distributions; this family
 contains all the commonly used discrete distributions.
 A complete classification of the Ord family of probability mass functions
 is related to the orthogonality of the corresponding Rodrigues polynomials.
 Also, for any random variable $X$ of this family and for
 any suitable function $g$ in $L^2(\RR,X)$, the article
 provides useful relationships between the Fourier
 coefficients of $g$ (with respect to the orthonormal
 polynomial system associated to $X$) and the Fourier
 coefficients of the forward difference of $g$ (with respect
 to another system of polynomials, orthonormal with respect
 to another distribution of the system). Finally,
 using these properties, a class of bounds for
 the variance of $g(X)$ is
 obtained, in terms of the forward differences of $g$. These
 bounds unify and improve several existing results.
\end{abstract}


\begin{classcode}
Primary 60E05, 	62E99, 05E35, 42A61; Secondary 60E15.
\end{classcode}

 \section{Indroduction}
 \label{section intro}

 \citet{Ord2} introduced the discrete analogue of Pearson's system.
 \citeauthor{Ord2}'s family contains all integer-valued random variables
 (rvs) whose probability mass function (pmf), $p$, satisfies
 \begin{equation}
 \label{eq Ord's Diff-eq}
 \frac{\Delta
 p(j-1)}{p(j)}=\frac{a-j}{(a+b_0)+(b_1-1)j+b_2j(j-1)}.
 \end{equation}
 Here, $\Delta$ is the forward difference operator
 and $p(j)$ is the pmf
 of the discrete rv
 $X$ and $j$ takes values in an integer
 interval. In the sequel, the term ``discrete rv'' is
 customized to mean ``integer-valued rv''.
 \Cref{eq Ord's Diff-eq}
 is the discrete analogue
 of Pearson's differential equation.
 \citeauthor{Ord2} classified these distributions
 according to the values of the
 parameters $a$, $b_0$, $b_1$ and $b_2$;
 see \cite[Table 2.1, p.~87]{JKB} .

 The present work is concerned with the {\it Cumulative Ord Family} of
 discrete distributions,
 defined as follows.
 \begin{definition}[Cumulative Ord Family]
 \label{def CO}
 Let $X$ be a discrete rv with finite mean $\mu$ and pmf $p(j)=\Pr(X=j)$,
 $j\in\ZZ$.
 We say that $X$ belongs to the Cumulative Ord family (or that $p$ belongs
 to
 the Cumulative Ord system)
 if there exists a quadratic $q(j)=\delta j^2+\beta j+\gamma$ such that
 \begin{equation}
 \label{eq CO}
 \sum_{k\le{j}}(\mu-k)p(k)=q(j)p(j)
 \quad \text{for all}  \
 j\in\ZZ.
 \end{equation}
 If \eqref{eq CO} is satisfied, we write $X\sim\COq$ or $p\sim\COq$,
 or more explicitly,
 $X$ or $p\sim\CO$.
 \end{definition}

 Let $X\sim\COq$. \citet{APP2} studied the orthogonal
 polynomials generated by a Rodrigues-type formula
 (see \Cref{theo orthogonal polynomials in CO} below)
 and based on these polynomials, they prove Stein-type covariance identities
 \cite[see][Eq.~(2.7), p.~512]{APP2}.
 First-order covariance identities, $k=1$, of that kind
 have been studied by \citet{SL2012} for the Ord, Katz
 as well as modified power series families of distributions.
 \citet{APP2}, using Bessel's inequality
 (for $m\in\NN=\{0,1,\ldots\}$ such that $\E|X|^{2m}<\infty$
 and $\E\left[q^{[m]}(X)|\Delta^{m}g(X)|\right]<\infty$),
 showed that
 \begin{equation}
 \label{eq Variance bounds Bessel-type}
 \Var[g(X)]\ge\sum_{k=1}^{m}
 \frac{\E^2\left[q^{[k]}(X)\Delta^{k}g(X)\right]}{k!\Pi_\delta^{[k]}(k-1) \E\left[q^{[k]}(X)\right]};
 \end{equation}
 the equality holds iff $g$ coincides with a polynomial of degree at most
 $m$ in the support of $X$.
 For $q^{[k]}$, $\Pi_\delta^{[k]}$ and $\Delta^{k}$,
 see \Cref{notations} below.
 Also, for $n\in\NN$ such that $\E|X|^{2n}<\infty$
 and $\E\left\{q^{[n]}(X)\left[\Delta^{n}g(X)\right]^2\right\}<\infty$,
 \citet{APP1}, by applying
 a discrete Mohr and Noll inequality,
 established some
 Poincar\'e-type
 bounds for the variance of $g(X)$, of the form
 \begin{equation}
 \label{eq existing Poincare}
 (-1)^{n}\left\{\Var\left[g(X)\right]-{S}_{n}\right\}\geq0,
 \quad \textrm{where} \
 {S}_n=\sum_{k=1}^{n}(-1)^{k-1}
 \frac{\E\left\{q^{[k]}(X)\left[\Delta^{k}g(X)\right]^2\right\}}{k!\Pi_\delta^{[k]}(0)};
 \end{equation}
 the equality holds iff $g$ is identified with a polynomial of
 degree at most $n$ in the support of $X$.

 We first present a simple example for illustrating
 the improvement achieved by the results of the present article.
 Let $X\sim\mathrm{Poisson}(\lambda)$. For $m=n=1$,
 \eqref{eq Variance bounds Bessel-type} and
 \eqref{eq existing Poincare}
 produce the
 double inequality
 \begin{equation}
 \label{eq existing bounds for Poisson m=n=1}
 \lambda\E^2[\Delta g(X)] \le \Var[g(X)] \le \lambda\E[\Delta g(X)]^2,
 \end{equation}
 where both equalities hold iff $g$ is a linear polynomial.
 Applying the results of the present paper
 (see \Cref{theo variance bounds}),
 we get the strengthened inequality
 \begin{equation}
 \label{eq new bound for Poisson m=n=1}
 \Var [g(X)] \le \frac{\lambda}{2}\E^2[\Delta g(X)]
 + \frac{\lambda}{2}\E[\Delta g(X)]^2,
 \end{equation}
 in which the equality holds iff $g$ is a polynomial of degree at most two.
 It is clear that the upper bound in \eqref{eq new bound for Poisson m=n=1}
 improves the upper bound in \eqref{eq existing bounds for Poisson m=n=1}
 and, in fact, it is strictly better,
 unless $g$ is linear.

 The rest of this paper is organized as follows.
 In \Cref{section preliminaries},
 we present some elementary properties of the
 Cumulative Ord (CO) family of distributions.
 In \Cref{section classification}, we give an algorithm
 (\Cref{algorithm}) which checks if a pair $(\mu;q)$ is admissible,
 according to \Cref{def (mu;q) admissible}.
 Also, we provide a detailed classification of the CO family.
 It turns out that, up to an (integer-valued) location
 transformation and/or multiplication by $-1$,
 there are six different types of pmfs,
 described in \Cref{table probabilities},
 while \Cref{section comparison} offers a comparison between \citeauthor{Ord1}'s
 Discrete Student distributions and that ones presented in this article.
 \Cref{section symmetric} presents the symmetric
 distributions  of the CO family of distributions;
 moreover, in this section, we define the
 noncentrality parameter as well as the degrees of freedom
 of a discrete Student distribution.
 In \Cref{section moment relations in CO},
 we show that for any $p\sim\COq$,
 the pmf $p_i\propto q^{[i]}p$ belongs to the CO system,
 under appropriate moment conditions.
 Also, using a known covariance identity, we obtain
 close-form expressions for $\Var(X_i)$ (where $X_i\sim p_i$) and
 for $\E\left[q^{[i]}(X)\right]$.
 Recurrence relations for the factorial moments are also given.
 In \Cref{section orthogonal polynomials in CO},
 we study the Rodrigues-type orthogonal polynomials of a rv that belongs to the CO Family.
 The main result of this section is that the forward differences of orthogonal polynomials
 of a pmf within the CO system are also orthogonal polynomials corresponding to
 another pmf of the system; see \Cref{lem Delta phi_k} and
 \Cref{theo Delta^m phi_k}.
 In \Cref{section L^2 completeness and expansions},
 we present expressions for the Fourier coefficients
 of a function $g$,
 with respect to the corresponding orthonormal polynomials.
 One of the important facts is that
 when the polynomials are dense in $L^2(\R,X)$,
 the expectation $\E\left\{q^{[n]}(X)[\Delta^ng(X)]^2\right\}$ can
 be expressed as a series (finite or infinite) in terms of
 the Fourier coefficients of $g$.
 This is utilized in \Cref{section applications to variance bounds}
 where, upon applying the
 series expansion for $\E\left\{q^{[n]}(X)[\Delta^ng(X)]^2\right\}$,
 we present a wide class of (upper/lower) bounds
 for $\Var[g(X)]$.

 \section{Preliminaries}
 \label{section preliminaries}
 In the present section, we investigate the basic
 properties of the CO family. In the sequel,
 we
 shall make use of the following notation.
 \begin{notation}
 \label{notations}
 \begin{enumerate}[itemsep=.5ex, wide, labelwidth=!, labelindent=0pt, label=\rm(\alph*), ref=\textcolor{black}{\alph*}]
 \item
 \label{notations(a)}
 $S\equiv S(X)=\{j\in\ZZ\colon\Pr(X=j)>0\}$ will denote the support of a
 discrete rv $X$. Also, we define
 $S_\circ\doteq S\smallsetminus\{\textrm{the lower endpoint of $S$}\}$ and
 $S^\circ\doteq S\smallsetminus\{\textrm{the upper endpoint of $S$}\}$.
 Of course, if $S$ does not have a finite lower (upper) endpoint,
 then $S_\circ=S$ ($S^\circ=S$);
 \item
 \label{notations(b)}
 For each real function $f$ and $k=0,1,\ldots$,
 we define $f^{[k]}(x)=f(x)\cdots f(x+k-1)$ and $f^{[-k]}(x)=1/[f(x-1)\cdots f(x-k)]$
 (provided that $f(x-1)\cdots f(x-k)\ne0$), with $f^{[0]}(x)=1$;
 \item
 \label{notations(c)}
 For $z\in\CC$ and $k=0,1,\ldots$,
 we define $[z]_k=z(z+1)\cdots(z+k-1)$ and $[z]_{-k}=[(z-1)\cdots(z-k)]^{-1}$
 (provided that $z\ne1,\ldots,k$), with $[z]_0=1$;
 \item
 \label{notations(d)}
 For $z\in\CC$ and $k\in\ZZ$,
 we define $(z)_k=(-1)^k[-z]_k$, provided that the quantity $[-z]_k$ is well-defined.
 That is, for $k=0,1,\ldots$, we have $(z)_k=z(z-1)\cdots(z-k+1)$
 and $(z)_{-k}=[(z+1)\cdots(z+k)]^{-1}$ (provided that $z\ne-1,\ldots,-k$), with $(z)_0=1$;
 \item
 \label{notations(e)}
 $\Delta^k=\Delta(\Delta^{k-1})$ with $\Delta^0=I$, the
 $k$th order forward difference operator;
 \item
 \label{notations(f)}
 For each $\delta\in\RR$ and $k\in\NN$, $\Pi_\delta^{[k]}(m)\doteq\prod_{j=m}^{k+m-1}(1-j\delta)$, $m\in\NN$, with $\Pi_\delta^{[0]}(m)=1$;
 \item
 \label{notations(g)}
 Let $\bm{z}=(z_1,z_2)$ and $\bm{w}=(w_1,w_2)\in\CC^2$. We denote $\bm{z}\eqcirc\bm{w}$ if $(z_1,z_2)=(w_1,w_2)$ or $(z_1,z_2)=(w_2,w_1)$.
 \end{enumerate}
 \end{notation}

 \begin{remark}
 \label{rem unimodal}
 For every non-degenerate discrete rv $X$ with finite mean $\mu$ and pmf $p$,
 the function $f(j)=\sum_{k\le{j}}(\mu-k)p(k)$, $j\in\ZZ$, is non-negative,
 unimodal (increases for $j\leq \lfloor\mu\rfloor$, the integer part of $\mu$, and then decreases)
 and takes its maximum value at the point $j=\lfloor\mu\rfloor$
 (if $\mu\in\ZZ$, the maximum value is attained at the points $\mu-1$ and $\mu$).
 \end{remark}

 \begin{definition}
 \label{def integer chain}
 A subset $S$ of $\ZZ$ is called an {\it integer chain} if for every $j_1,j_2\in S$
 with $j_1\le j_2$, we have $[j_1,j_2]\cap\ZZ\subseteq S$.
 \end{definition}

 \begin{remark}
 \label{rem integer chain}
 An integer chain $S$ of $\ZZ$ can be written as $\{\alpha,\ldots,\omega\}$,
 where $\alpha\in\ZZ\cup\{-\infty\}$ and $\omega\in\ZZ\cup\{\infty\}$,
 in the sense that
 $\{\alpha,\ldots,\infty\}=\{\alpha,\alpha+1,\ldots\}$ when $\alpha>-\infty$,
 $\{-\infty,\ldots,\omega\}=\{\ldots,\omega-1,\omega\}$ when $\omega<\infty$
 and $\{-\infty,\ldots,\infty\}=\ZZ$.
 \end{remark}

 From \Cref{def CO}, we can prove the following lemma.

 \begin{lemma}
 \label{lem S:q:S_o:S^o}
 Let $X\sim\COq$. Then:
 \begin{enumerate}[itemsep=.5ex, wide, labelwidth=!, labelindent=0pt, label=\rm(\alph*), ref=\textcolor{black}{\alph*}]
 \item
 \label{lem S:q:S_o:S^o(a)}
 The rv $X$ is supported on an integer chain, denoted by
 $S=\{\alpha,\ldots,\omega\}$,
 where $\alpha\in\ZZ\cup\{-\infty\}$ and $\omega\in\ZZ\cup\{\infty\}$
 with $\alpha\le\omega$. Thus,
 $S_\circ=\{\alpha+1,\ldots,\omega\}$ and $S^\circ=\{\alpha,\ldots,\omega-1\}$.
 Note that if $\alpha=\omega$, then $S_\circ=S^\circ=\varnothing$;
 \item
 \label{lem S:q:S_o:S^o(b)}
 $q(j)>0$ for all $j\in S^\circ$;
 \item
 \label{lem S:q:S_o:S^o(c)}
 If $\omega<\infty$, then $q(\omega)=0$. If $\alpha=0$, then
 $\mu=q(0)=\gamma$ and, in general, if $\alpha>-\infty$,
 then $q(\alpha)=\mu-\alpha$;
 \item
 \label{lem S:q:S_o:S^o(d)}
 For every $r\in\ZZ$, the rv $Y=X+r$ follows $\mathrm{CO}(\mu_Y=\mu+r;q_Y(j)=q(j-r))$;
 \item
 \label{lem S:q:S_o:S^o(e)}
 The rv $W=-X$ follows
 $\mathrm{CO}(\mu_W=-\mu;q_W(j)=q(-j)-j-\mu)$;
 \item
 \label{lem S:q:S_o:S^o(f)}
 $\underline{q}(j)>0$ for all $j\in S_\circ$, where $\underline{q}(j)\doteq q(j)+j-\mu$;
 \item
 \label{lem S:q:S_o:S^o(g)}
 $p(j)=r(j-1)p(j-1)$ for all $j\in S_\circ$, where $r(j)\doteq q(j)\left/\underline{q}(j+1)\right.$, $j\in S^\circ$.
 \end{enumerate}
 \end{lemma}
 \begin{proof}
 \eqref{lem S:q:S_o:S^o(a)}, \eqref{lem S:q:S_o:S^o(b)} and \eqref{lem S:q:S_o:S^o(c)} are obvious by \eqref{eq CO} and \Cref{rem unimodal}.

 \eqref{lem S:q:S_o:S^o(d)}
 The rv $Y$ has mean $\mu_Y=\mu+r$ and
 support $S(Y)=r+S=\{r+j\colon j\in S\}$ (integer chain).
 Observe that
 $\sum_{k\le{j}}(\mu_Y-k)\Pr(Y=k)=\sum_{k\le{j}}[\mu-(k-r)]\Pr(X=k-r)=\sum_{s\le{j-r}}(\mu-s)\Pr(X=s)=q(j-r)\Pr(X=j-r)=q(j-r)\Pr(Y=j)$.

 \eqref{lem S:q:S_o:S^o(e)}
 The rv $W$ has mean $\mu_W=-\mu$ and support $S(W)=-S=\{-j\colon j\in S\}$. Now, write
 $\sum_{k\le{j}}(\mu_W-k)\Pr(W=k)=\sum_{k\le{j}}(-\mu-k)\Pr(X=-k)=-\sum_{k\le{j}}[\mu-(-k)]\Pr(X=-k)=-\sum_{s\ge-j}(\mu-s)\Pr(X=s)=(-\mu-j)\Pr(X=-j)-\sum_{s>-j}(\mu-s)\Pr(X=s)$
 and observe that $\sum_{s\le-j}(\mu-s)\Pr(X=s)+\sum_{s>-j}(\mu-s)\Pr(X=s)=0$. Thus,
 $\sum_{k\le{j}}(\mu_W-k)\Pr(W=k)=(-\mu-j)\Pr(X=-j)+\sum_{s\le-j}(\mu-s)\Pr(X=s)=(-\mu-j)\Pr(X=-j)+q(-j)\Pr(X=-j)=[q(-j)-j-\mu]\Pr(W=j)$.

 \eqref{lem S:q:S_o:S^o(f)}
 It is obvious that $j\in S_\circ(X)\Leftrightarrow-j\in S^\circ(-X)$.
 Thus, from \eqref{lem S:q:S_o:S^o(b)} and \eqref{lem S:q:S_o:S^o(e)}, we see that for every $j\in S_\circ$,  $q_{-X}(-j)>0$, i.e.,
 $q(j)+j-\mu>0$.

 \eqref{lem S:q:S_o:S^o(g)}
 From \eqref{lem S:q:S_o:S^o(f)}, it follows that $r(j)$ is well-defined. For each $j\in S_\circ$, \eqref{eq CO} gives $q(j)p(j)=\sum_{k\le{j}}(\mu-k)p(k)=q(j-1)p(j-1)+(\mu-j)p(j)$, that is, $\underline{q}(j)p(j)=q(j-1)p(j-1)$.
 \end{proof}

 Now, we present a useful lemma concerning the existence of moments
 (see \cite[p.~176]{APP1}).

 \begin{lemma}
 \label{lem moments and delta}
 Let $X\sim\CO$. If $S(X)$ is finite or $\delta\le0$, then $X$ has finite
 moments of any order.
 Furthermore, if $S(X)$ is infinite and $\delta>0$, then
 $X$ has finite moments of any order $\theta\in[0,1+1/\delta)$,
 while $\E|X|^{1+1/\delta}=\infty$.
 \end{lemma}

 \begin{remark}
 \label{rem finite S: delta>0}
 We can find a rv $X\sim\CO$ with $\delta>0$ and finite support $S$
 (with cardinality $|S|\ge3$).
 However, the inequality $\delta<[2(|S|-2)]^{-1}$ should be
 necessarily satisfied in this
 case; see \Cref{sssection delta>0: finite S}.
 \end{remark}

 \begin{lemma}
 \label{lem lim_j->oo q(j)p(j)}
 Suppose a discrete rv $X$ is supported on an infinite integer chain $S(X)$,
 has finite mean $\mu$, and satisfies the relation $\sum_{k\le{j}}(c-k)p(k)=q(j)p(j)$,
 $j\in S(X)$, where $c$ is a constant and $q$ is a polynomial of degree at most two.
 If $S(X)$ is upper (resp., lower) unbounded, then $q(j)p(j)\to0$ as $j\to \infty$
 (resp., $j\to-\infty$).
 \end{lemma}
 \begin{proof}
 If $q$ is a linear polynomial, then the result
 is obvious because $\E|X|<\infty$. We shall
 examine only the case where
 $q(j)=\delta j^2+\beta j+\gamma$ with $\delta\ne0$.
 If $S(X)$ is upper unbounded, the quantity
 $\sum_{k\in\ZZ}(c-k)p(k)=\lim_{j\to\infty}\sum_{k\le{j}}(c-k)p(k)$
 is well-defined, since $\E|X|<\infty$.
 Thus, $\lim_{j\to\infty}p(j)q(j)=C\in\RR$.
 Observe that
 $\lim_{j\to\infty}j^2p(j)=\lim_{j\to\infty}\left[\left.j^2\right/q(j)\right]q(j)p(j)=C/\delta\ge0$.
 Assuming $C>0$, we can find an integer $j_0>0$ such that $j p(j)>C/(2\delta j)$
 for all $j\ge j_0$.
 Thus,
 $\E|X|\ge\sum_{j\ge{j_0}} j p(j)\ge \sum_{j\ge{j_0}}C/(2\delta j)=\infty$,
 a contradiction.
 For the lower unbounded case, we can use analogous arguments.
 \end{proof}

 The result of \Cref{lem lim_j->oo q(j)p(j)}
 applies to all rvs of the CO family whose
 support is infinite.
 For this family, the results of the above lemma can be generalized;
 see \Cref{prop lim_j->oo j^2ip(j)}.

 \begin{remark}
 \label{rem c=mu for upper unbouded S}
 In \Cref{lem lim_j->oo q(j)p(j)},
 if the support does not have a finite upper bound, then
 the constant $c$ is necessarily the mean $\mu$ of $X$, i.e., the rv $X$ belongs to
 the CO family. However,
 if the support has a finite upper bound,
 it may happen that
 $c\ne\mu$. For example,
 let $X\sim\mathrm{Poisson}(\lambda)=\mathrm{CO}(\mu=\lambda;q(j)=\lambda)$.
 Then, from \Cref{lem S:q:S_o:S^o}\eqref{lem S:q:S_o:S^o(e)},
 the rv $W=-X\sim\mathrm{CO}(\mu_W=-\lambda;q_W(j)=-j)$,
 i.e., $\sum_{k=-\infty}^{j}(\mu_W-k)p_W(k)=q_W(j)p_W(j)$ for all $j\in\ZZ$.
 Now, consider the truncated of $W$, say $V$, at the point zero
 [$p_{V}(j)=\vartheta p_W(j)$, $j=-1,-2,\ldots$].
 Then, for each
 $j\in S(V)=\{-1,-2,\ldots\}$,
 we have $\sum_{k=-\infty}^{j}(\mu_W-k)p_{V}(k)=\vartheta\sum_{k=-\infty}^{j}(\mu_W-k)p_W(k)
 =\vartheta q_W(j)p_W(j)=(-j)p_{V}(j)$
 and $\E(V)\ne\mu_W$
 (note that this $V$ does not belong
 to the CO family).
 \end{remark}

 Now, we compare the CO system,
 i.e., the pmfs satisfying \cref{eq CO},
 with the ordinary \citeauthor{Ord2} system, i.e.,
 the pmfs satisfying the \citeauthor{Ord2}'s difference
 \cref{eq Ord's Diff-eq}.

 \begin{proposition}
 \label{prop CO against Ord system}
 Assume that a discrete rv $X$ has pmf $p$ and finite mean.
 Then, \eqref{prop CO against Ord system(a)} and \eqref{prop CO against Ord system(b)} are equivalent, where
 \begin{enumerate}[nolistsep, leftmargin=19pt, label=\rm(\alph*), ref=\textcolor{black}{\alph*}]
 \item
 \label{prop CO against Ord system(a)}
     $X\sim\COq$;

 \item
 \label{prop CO against Ord system(b)}
     \begin{enumerate}[nolistsep, leftmargin=19pt, label=\rm(\roman*), ref=\textcolor{black}{\roman*}]
     \item
     \label{prop CO against Ord system(b)i}
     The support $S$ of $X$ is an integer chain,
     $S=\{\alpha,\ldots,\omega\}$ with $\alpha\leq \omega$,
     $\alpha\in\ZZ\cup\{-\infty\}$, $\omega\in\ZZ\cup\{\infty\}$;
     \item
     \label{prop CO against Ord system(b)ii}
     There exist polynomials $p_1$ (of degree at most one) and
     $q$ (of degree at most two) such that $[\Delta p(j-1)]q(j-1)=p_1(j)p(j)$, for all $j\in S$;
     \item
     \label{prop CO against Ord system(b)iii}
     For the above polynomials, there exists a constant $\mu$ such that
     $p_1(j)+\Delta q(j-1)=\mu-j$, for all $j\in S$.
     If $\omega<\infty$, then it is further required that
     $\mu=\E(X)$.
     \end{enumerate}
 \end{enumerate}
 \end{proposition}

 \begin{proof}
It is obvious that \eqref{prop CO against Ord system(a)} implies \eqref{prop CO against Ord system(b)}.
We now prove that \eqref{prop CO against Ord system(b)} implies \eqref{prop CO against Ord system(a)}.
Note that we do not assume that $\mu$ is the mean
 of $X$. Combining \eqref{prop CO against Ord system(b)ii} and
 \eqref{prop CO against Ord system(b)iii}, we have, after some algebra, that
 $p(j)(\mu-j)=\Delta[p(j-1)q(j-1)]$. Fix an integer $i$ with $i\le j$.
 Then, $\sum_{k=i}^{j}p(k)(\mu-k)=p(j)q(j)-p(i-1)q(i-1)$.
 If $\alpha>-\infty$, we choose $i=\alpha$ and since $p(\alpha-1)=0$, we have
 $\sum_{k\le{j}}p(k)(\mu-k)=\sum_{k=\alpha}^{j}p(k)(\mu-k)=p(j)q(j)$,
 for all $j\in S$.
 For $\alpha=-\infty$, since $\E|X|<\infty$, the quantity
 $\sum_{k\le{j}}(\mu-k)p(k)=\lim_{i\to-\infty}\sum_{k=i}^{j}(\mu-k)p(k)$
 is well-defined.
 Thus, $\lim_{i\to-\infty}p(i-1)q(i-1)=C\in\RR$.
 If $q$ is of degree at most one, then it is obvious that $C=0$.
 If $q(j)=\delta j^2+\beta j+\gamma$, $\delta\ne0$,
 we use the same arguments as in the proof of \Cref{lem lim_j->oo q(j)p(j)}
 and we conclude that $C=0$.
 Thus, in any case, $\sum_{k\le{j}}(\mu-k)p(k)=q(j)p(j)$, $j\in S$.
 \Cref{lem lim_j->oo q(j)p(j)} and
 \Cref{rem c=mu for upper unbouded S} then complete the proof.
 \end{proof}

 \begin{remark}
 \label{rem CO against Ord system}
 All assumptions of \Cref{prop CO against Ord system}\eqref{prop CO against Ord system(b)}
 are necessary for a rv to lie in the CO family:
 It is obvious that \eqref{prop CO against Ord system(b)ii} is necessary for all $j\in S$.
 Regarding the assumption \eqref{prop CO against Ord system(b)i}, consider the rv $X$ with pmf $p(0)=p(2)=0.5$
 and observe that
 \eqref{prop CO against Ord system(b)ii}
 and \eqref{prop CO against Ord system(b)iii} are fulfilled for
 $p_1(j)=2-j$, $q(j)=1-j$ and $\mu=1=\E(X)$. But, this rv
 does not belong to the CO family.
 Regarding the assumption \eqref{prop CO against Ord system(b)iii}, consider the truncated Poisson
 $X\sim p(j)\propto \left.\lambda^j\right/j!$, $j=0,1,\ldots,N$.
 Observe that \eqref{prop CO against Ord system(b)i}--\eqref{prop CO against Ord system(b)iii}
 are satisfied for $p_1(j)=\lambda-j$,
 $q(j)=\lambda$ and $\mu=\lambda$.
 However, since $\lambda\ne\E(X)$, this rv does not belong
 to the CO family;
 this is so because $S$ has a finite
 upper endpoint.
 \end{remark}

 \section{A complete classification of the Cumulative Ord family}
 \label{section classification}
 In this section, we classify the distributions of the CO family. The
 classification is based on the mean $\mu$ and the
 parameters of the quadratic $q$. The most important role
 is played by the parameter $\delta$, the coefficient of the square
 power of $q$.

 The natural question is to ask whether the mean $\mu$
 and the quadratic $q$, together, characterize the
 distribution. The answer is given in the following proposition.

 \begin{proposition}
 \label{prop the support is defined by mu and q}
 Suppose the rv $X$ follows the $\COq$ distribution. Then:
 \begin{enumerate}[nolistsep, wide, labelwidth=!, labelindent=0pt, label=\rm(\alph*), ref=\textcolor{black}{\alph*}]
 \item
 \label{prop the support is defined by mu and q(a)}
 The support $S(X)$ is uniquely determined by $\mu$ and $q$;
 \item
 \label{prop the support is defined by mu and q(b)}
 The distribution is characterized by the pair $(\mu;q)$.
 \end{enumerate}
 \end{proposition}

 \begin{proof}
 \eqref{prop the support is defined by mu and q(a)}
 First, we consider the special case when $\mu\in\ZZ$ and
 $q(\mu)=0$. From \Cref{rem unimodal}, it follows easily
 that the rv takes the value $\mu$ with probability 1.
 Otherwise, we define $N\equiv N(\mu;q)=\{\textrm{the first
 $j\in \ZZ \cap[\mu,\infty)$ such that $q(N)\le0$}\}$. If
 $N=\infty$, then the support $S(X)$ does not have a finite upper
 endpoint, otherwise
 the value $N$ is the upper endpoint of $S(X)$ and, then,
 $q(N)=0$ (otherwise the pair $(\mu;q)$ could not satisfy
 the relation \eqref{eq CO}).
 Regarding the lower endpoint: The rv $-X$ follows
 $\mathrm{CO}(\mu_{-X};q_{-X})$ with $S(-X)=-S(X)$,
 $\mu_{-X}=-\mu_{X}$ and $q_{-X}(j)=q(-j)-j-\mu$. As before,
 we can determine  $N'\equiv N(\mu_{-X};q_{-X})$.
 If $N'=\infty$, then the
 support $S(-X)$ does not have a finite upper endpoint, i.e., the support
 $S(X)$ does not have a finite lower endpoint. Otherwise, the value $N'$ is the
 upper endpoint of $S(-X)$, i.e., $-N'$ is the
 lower endpoint of $S(X)$.

 \eqref{prop the support is defined by mu and q(b)}
 Consider the pmf $p$ of $X$ and its support, $S(X)$,
 which is determined by the pair $(\mu;q)$.
 Now, let $\tilde{p}\sim\COq$.
 Consider the function $r$ of
 \Cref{lem S:q:S_o:S^o}\eqref{lem S:q:S_o:S^o(g)} and fix $j_0\in S(X)$.
 For every $k\in\ZZ$ with
 $j_0+k\in S(X)$, it follows that
 $p(j_0+k)=r^{[k]}(j_0)p(j_0)$
 and $\tilde{p}(j_0+k)=r^{[k]}(j_0)\tilde{p}(j_0)$.
 Consequently, $\tilde{p}\propto p$, and so
 $\tilde{p}=p$ because $p$ and $\tilde{p}$ are pmfs.
 \end{proof}

 \begin{definition}
 \label{def (mu;q) admissible}
 Let $\mu\in\RR$ and $q(j)=\delta j^2 + \beta j +\gamma$.
 We say that the pair $(\mu;q)$ is {\it admissible} if there exists a
 pmf $p$ in the CO system such that $p\sim\COq$.
 \end{definition}

 Now, the natural question is ``How one can check the
 admissibility of a given pair $(\mu;q)$?'' Also, if a pair
 is admissible, how can the corresponding pmf be
 obtained by this pair? The answer is given in \Cref{algorithm}.

\begin{algorithm}
\begin{minipage}[t]{1\columnwidth}
\small
\caption{Admissibility of the pair $(\mu;q)$.}
 \label{algorithm}
\begin{algorithmic}[1]
\State
 Consider the polynomial $\underline{q}(j)$ of \Cref{lem S:q:S_o:S^o}\eqref{lem S:q:S_o:S^o(f)} and define
 \[
 \alpha\doteq\sup\left\{j\in(-\infty,\mu]\cap\ZZ\colon \underline{q}(j)=0\right\},
 \quad
 \omega\doteq\inf\{j\in[\mu,\infty)\cap\ZZ\colon q(j)=0\},
 \]
 noting that $\sup\{\varnothing\}=-\infty$ and
 $\inf\{\varnothing\}=\infty$. Next, define
$S=[\alpha,\omega]\cap\ZZ$. The pair $(\mu;q)$
is admissible if and only if $q(j)>0$ for all $j\in
S^\circ$ and $\underline{q}(j)>0$ for all $j\in S_\circ$.
If $(\mu;q)$ is admissible, go to Step 2, end otherwise.

\State
Let $p$ be the corresponding pmf to the pair $(\mu;q)$ in the CO system. By application of \Cref{lem S:q:S_o:S^o}\eqref{lem S:q:S_o:S^o(g)},
we get\footnotemark:
 \begin{subequations}
 \label{eq pmf a:omega}
 \begin{align}
 &\textrm{if } \alpha>-\infty,
 &&p(\alpha+i)\propto r^{[i]}(\alpha),\quad\hspace{1.75ex} i=0,1,\ldots,\omega-\alpha;\qquad\label{subeq pmf a>-infty}\\
 &\textrm{if }  \alpha=-\infty \textrm{ and } \omega<\infty,\qquad\qquad
 &&p(\omega-i)\propto r^{[-i]}(\omega),\quad i=0,1,\ldots; \label{subeq pmf a=-infty: omega<infty}\\
 &\textrm{if }  \alpha=-\infty \textrm{ and } \omega=\infty,
 &&p(j)\propto r^{[j]}(0),\qquad\quad\hspace{1ex} j\in\ZZ. \label{subeq pmf a=-infty: omega=infty}
 \end{align}
  \end{subequations}
\end{algorithmic}
\mpfootnotes[1]
\footnote{For the pmfs in \eqref{eq pmf a:omega}, we observe the following.
 Since \Cref{lem S:q:S_o:S^o}\eqref{lem S:q:S_o:S^o(g)} holds, as in the analysis presented in \cite[Lemma 4.1, pp.~176--178]{APP1},
 one can see that $\sum_{j\in S}|j|p(j)<\infty$ (so the mean of $X$ is finite), and hence, $\sum_{j\in S}p(j)<\infty$.
 By construction, these pmfs satisfy the relation
 $[\Delta p(j-1)]q(j-1)=[-\Delta q(j-1)-j+\mu]p(j)$, $j\in S$.
 As in the proof of \Cref{prop CO against Ord system},
 we get $\sum_{k\le{j}}(\mu-k)p(k)=q(j)p(j)$ for all $j\in S$.
 If $\omega=\infty$, then $q(j)p(j)\to0$ as $j\to\infty$; see
 \Cref{lem lim_j->oo q(j)p(j)}; thus,
 $\sum_{k\in{S}}(\mu-k)p(k)=0$. If $\omega<\infty$, then $q(\omega)=0$, see Step 1 of \Cref{algorithm}, and so $\sum_{k\in{S}}(\mu-k)p(k)=\sum_{k\le{\omega}}(\mu-k)p(k)=q(\omega)p(\omega)=0$.
 In both cases, $\omega<\infty$ and $\omega=\infty$, $\mu$ is the mean value.}
\end{minipage}
\end{algorithm}

 Next, we present a detailed classification of the CO system.

 \subsection[The case $\delta=0$]{The case $\bm{\delta=0}$}
 \label{ssection delta=0}
 We have to further distinguish between the cases $\beta=0$ and $\beta\ne0$.

 \subsubsection{The subcase $\beta=0$ {\rm(Poisson-type distributions)}}
 \label{sssection delta=0:beta=0}
 Then, $q(j)=\gamma>0$, $j\in S$. The support $S$ does not have a finite
 upper endpoint, but it must have a finite lower endpoint (because the
 quadratic $q_{-X}(j)=\gamma-j-\mu$ of $-X$, see
 \Cref{lem S:q:S_o:S^o}\eqref{lem S:q:S_o:S^o(e)}, takes negative values for large
 values of $j$). Without lost of generality, we assume
 $S=\NN$. Since $\alpha=0$, $\mu=\gamma$. Observe that the
 Poisson distribution with parameter $\lambda=\gamma$
 follows $\mathrm{CO}(\gamma;0,0,\gamma)$.
 Using \Cref{prop the support is defined by mu and q},
 we have $X\sim\mathrm{Poisson}(\gamma)$.

 \subsubsection{The subcase $\beta\ne0$}
 \label{sssection delta=0:beta neq 0}
 We have the following sub-subcases.

 \paragraph{The sub-subcase $\beta>0$ {\rm(Negative Binomial-type distributions)}}
 \label{beta-ne0(i)}
 The support $S$ does not have a finite upper endpoint, but it has
 a lower one. Again, we may assume that $S=\NN$ (of course
 $\mu>0$). Since $\alpha=0$, $q(j)=\beta j+\mu$. Consider
 the Negative Binomial distribution with parameters
 $r=\beta/\mu>0$ and $p=1/(1+\beta)\in(0,1)$, i.e.,
      \[
      p(j)=\frac{[r]_j}{j!}p^r(1-p)^j,\quad j=0,1,\ldots,
      \]
 which follows $\mathrm{CO}(\mu;0,\beta,\mu)$.
 From \Cref{prop the support is defined by mu and q},
 $X\sim\mathrm{NB}(r=\beta/\mu,p=1/(1+\beta))$.

 \paragraph{The sub-subcase $-1<\beta<0$ {\rm(Binomial-type distributions)}}
 \label{beta-ne0(ii)}
 The support $S$ has a finite upper endpoint. Also,
 $q_{-X}(j)=-(1+\beta)j+\gamma-\mu$ where, $-(1+\beta)<0$.
 Thus, $S$ has a finite lower endpoint. Assume that
 $S=\{0,1,\ldots,N\}$ and $0<\mu<N$. Since $q(0)=\mu$,
 $q(N)=0$ and $q$ is a linear polynomial, we get
 $q(j)=\mu(N-j)/N$. Consider the Binomial
 distribution with parameters $N$ and $p=\mu/N$, i.e.,
      \[
      p(j)={N\choose{j}}p^j(1-p)^{N-j},\quad j=0,1,\ldots,N,
      \]
      which follows $\mathrm{CO}(p;0,-p,Np)$.
      From \Cref{prop the support is defined by mu and q},
      we see that $X\sim \mathrm{Bin}(N,p=\mu/N)$.

 \paragraph{The sub-subcase $\beta=-1$ {\rm(Poisson-type distributions)}}
 \label{beta-ne0(iii)}
 Here, $q_{-X}(j)=\gamma-\mu$ (constant). Thus, this
 sub-subcase is the negative of the case \ref{sssection
 delta=0:beta=0}.

 \paragraph{The sub-subcase $\beta<-1$ {\rm(Negative Binomial-type distributions)}}
 \label{beta-ne0(iv)}
 In this case, $q_{-X}(j)=-(1+\beta)j+\gamma-\mu$, $-(1+\beta)>0$. This case is the negative of the
 case \ref{beta-ne0(i)}.

 \subsection[The case $\delta<0$]{The case $\bm{\delta<0}$ \emph{(Negative Hyper\-geometric-type distributions)}}
 \label{ssection delta<0}
 It is obvious that $S$ is finite. Without loss
 of generality, assume that $S=\{0,1,\ldots,N\}$ with
 $0<\mu<N$. Since $q(N)=0$ and $q(0)=\mu$, it follows that
 \mbox{$q(j)=\delta[\mu/(N\delta)-j](N-j)$.}
 Consider now the
 Negative Hypergeometric distribution with parameters
 $N\in\{1,2,3,\ldots\}$, $r=-\mu/(N\delta)>0$ and
 $s=(\mu-N)/(N\delta)>0$, i.e., with pmf
 \[
 p(j)={N\choose{j}}\frac{(-r)_j(-s)_{N-j}}{(-r-s)_N},\quad j=0,1,\ldots,N.
 \]
 This follows $\mathrm{CO}\left(\frac{Nr}{r+s};\frac{-1}{r+s},\frac{N-r}{r+s},\frac{Nr}{r+s}\right)$.
 From \Cref{prop the support is defined by mu and
 q}, $X\sim\mathrm{NHgeo}(N;r,s)$.

 \subsection[The case $\delta>0$]{The case $\bm{\delta>0}$}
\label{ssection delta>0}
 We study the following subcases, relating to the support.

 \subsubsection{Finite $S$ {\rm(Hyper\-geometric-type distributions)}}
 \label{sssection delta>0: finite S}
Set $S=\{0,1,\ldots,N\}$ and $0<\mu<N$.
As in \Cref{ssection delta<0} ($\delta<0$), it follows that
$q(j)=\delta[\mu/(N\delta)-j](N-j)$. From \Cref{lem
S:q:S_o:S^o}\eqref{lem S:q:S_o:S^o(b)}, we get $\mu/(N\delta)-(N-1)>0$,
or equivalently $\delta<{\mu}/[N(N-1)]$. Now, since
$\underline{q}(j)=j[\delta{j}+(1-\mu/N-N\delta)]$, from
\Cref{lem S:q:S_o:S^o}\eqref{lem S:q:S_o:S^o(f)}, we have that
$\delta+(1-\mu/N-N\delta)>0$, or equivalently
$\delta<(N-\mu)/[N(N-1)]$. So,
 \[
 0<\delta<{\min\{\mu,N-\mu\}}/[N(N-1)].
 \]
Note that
${\min\{\mu,N-\mu\}}/[N(N-1)]\le1/[2(N-1)]$.
Consider the Hypergeometric distribution with parameters
$N\in\{2,3,\ldots\}$, $r=\mu/(N\delta)>N-1$ and
$s=(N-\mu)/(N\delta)>N-1$, with pmf
 \[
 p(j)={N\choose{j}}\frac{(r)_j(s)_{N-j}}{(r+s)_N},\quad j=0,1,\ldots,N.
 \]
Thus,  $p\sim\mathrm{CO}\left(\frac{Nr}{r+s};\frac{1}{r+s},-\frac{r+N}{r+s},\frac{Nr}{r+s}\right)$.
From \Cref{prop the support is defined by mu and
q}, $X\sim\mathrm{Hgeo}(N;r,s)$.

 \subsubsection{One-side infinite $S$ {\rm(Generalized Inverse Polya or Discrete $F$-type distributions)}}
 \label{sssection delta>0: one-sided infinite S}
 First, we give an example. Let $q(j)=j^2+1$ and $\mu=1$. It
 follows that $\underline{q}(j)=j(j+1)$. Step 1 of \Cref{algorithm} gives
 $\alpha=0$ and $\omega=\infty$, namely, $S=\NN$, and the
 pair $(\mu;q)=(1;j^2+1)$ is admissible. Using \eqref{subeq pmf a>-infty}, we find
 the pmf $p\sim\mathrm{CO}(1;1,0,1)$ which is
 $p(j)=[{\pi}/{\sinh(\pi)}]\left.\left[\prod_{k=0}^{j-1}(k^2+1)\right]\right/[j!(j+1)!]$,
 $j\in\NN$; this pmf can be written as
 \begin{equation}
 \label{eq pmf of CO(1;1:0:1)}
 p(j)=\frac{\pi}{\sinh(\pi)}\frac{[\imath]_j[-\imath]_j}{j!(j+1)!},\quad
 j=0,1,\ldots,
 \end{equation}
 where $\imath$ is the complex unity.

 Let us consider the general case when $S=\NN$ and $\mu>0$.
 Since $\alpha=0$, the quadratic $q$ is of the form $q(j)=\delta
 j^2+\beta j+\mu$. Write $q(j)=\delta(j+z_1)(j+z_2)$, where
 $-z_1,-z_2$ are the complex roots of $q$. Since $q(j)>0$
 for every $j\in S$, we get
 $(z_1,z_2)\in\CCC_2\subset\CC^2$, where
 \[
 \textstyle
 \CCC_2\doteq\{(z,\bar{z})\colon
 z\in\CC\smallsetminus\RR\}\cup(0,\infty)^2\cup
 \left\{\bigcup_{n=0}^{\infty}(-n-1,-n)^2\right\}.
 \]
 Observe that $\underline{q}(j)=\delta j(j+\rho)$, where $\rho=(\delta+\beta+1)/{\delta}$.
 It is also required that $\underline{q}(j)>0$ for all $j\in\NN^*=\{1,2,\ldots\}$,
 that is, $\rho>-1$, or equivalently, $\beta>-\delta-1$.
 Under the above restrictions, the pair $(\mu;q)$ is admissible. Step 2 of
 the algorithm then yields
 \begin{equation}
 \label{eq pmf of inv Polya}
 p(j)=\frac{\varGamma(\rho-z_1)
 \varGamma(\rho-z_2)}{\varGamma(\rho)\varGamma(\rho')}
 \frac{[z_1]_j[z_2]_j}{j![\rho]_j},\quad j=0,1,\ldots,
 \end{equation}
 where $\rho'=1+1/\delta$.
 The substitution $z_{1,2}\mapsto\pm\imath$
 and $\rho\mapsto2$
 in \eqref{eq pmf of inv Polya} yields
 \eqref{eq pmf of CO(1;1:0:1)}.

 \subsubsection{Two-side infinite $S$ {\rm(Discrete Student-type distributions)}}
 \label{sssection delta>0: two-sided infinite S}
 First, we give an example. Let $q(j)=j^2+1$ and $\mu=0$. It
 follows that $\underline{q}(j)=j^2+j+1$. Applying \Cref{algorithm}, Step 1 gives
 $\alpha=-\infty$ and $\omega=\infty$, namely, $S=\ZZ$, and
 the pair $(\mu;q)=(0;j^2+1)$ is admissible.
 Applying \eqref{subeq pmf a=-infty: omega=infty}, the pmf of $\mathrm{CO}(0;1,0,1)$ distribution is
 \begin{equation}
 \label{eq pmf of CO(0;1:0:1)}
 p(j)\propto\frac{[\imath]_j[-\imath]_j}
 {\left[{3}/{2}+\imath{\surd{3}}/{2}\right]_j
 \left[{3}/{2}-\imath{\surd{3}}/{2}\right]_j},
 \quad j\in\ZZ.
 \end{equation}
 Note that the above choice of $(\mu;q)$ forces $S$ to be
 the entire $\ZZ$.

 For the general case, let $\mu\in\RR$,
 $q(j)=\delta(z_1+j)(z_2+j)$ and
 $\underline{q}(j)=\delta(w_1+j)(w_2+j)$, where $-z_1$,
 $-z_2$ and $-w_1$, $-w_2$ are the complex roots of $q$ and
 $\underline{q}$, respectively.
 Writing $\bm{z}=(z_1,z_2)$ and $\bm{w}=(w_1,w_2)$,
 since $q(j)>0$ and
 $\underline{q}(j)>0$ for all $j\in\ZZ$, it follows that
 $\bm{z},\bm{w}\in\tilde{\CCC}_2$, where
 \[
 \tilde{\CCC}_2\doteq\{(z,\bar{z})\colon
 z\in\CC\smallsetminus\RR\}\cup\left\{\mbox{$\bigcup$}_{n\in\ZZ}(n,n+1)^2\right\}.
 \]
 Note that the pair $(w_1,w_2)$ is a function of
 $(\mu;q)$, i.e., a function of the values $\mu$, $\delta$,
 $z_1$ and $z_2$. The pair $(\mu;q)$ is admissible; see Step
 1 of \Cref{algorithm}. From Step 2, we
 obtain a formula for the pmf as follows (cf.\ \cite{Ord1}):
 \begin{equation}
 \label{eq pmf delta>0: S=Z}
 p(j)\propto\frac{[z_1]_j[z_2]_j}{[w_1+1]_j[w_2+1]_j},\quad j\in\ZZ.
 \end{equation}
 Substituting $z_{1,2}\mapsto\pm\imath$ and
 $w_{1,2}\mapsto1/2\pm\imath\surd{3}/2$
 in \eqref{eq pmf delta>0: S=Z}, we
 obtain \eqref{eq pmf of CO(0;1:0:1)}.

 All the above possibilities (\Cref{ssection delta=0,ssection delta<0,ssection delta>0})
 are summarized in \Cref{table probabilities} below.

  \begin{remark}
 \label{rem |S|=2}
 It is easy to check that if the cardinality $|S|$ of the
 support equals $2$, then different types lead to
 identical distributions (since every such rv
 is Bernoulli).
 \end{remark}

 \begin{sidewaystable}
 \footnotesize \caption[Probabilities of the cumulative Ord
 family]{Probabilities of the cumulative Ord family}
 \label{table probabilities}
 \footnotesize
 \begin{threeparttable}
 \begin{tabular*}{\textwidth}
 {@{\hspace{0ex}}@{\extracolsep{\fill}}c@{\hspace{0ex}}c@{\hspace{0ex}}c@{\hspace{0ex}}c@{\hspace{0ex}}c@{\hspace{0ex}}c@{\hspace{0ex}}r@{\hspace{0ex}}}
 \addlinespace
 \toprule
 \XA{\bf type}{\bf notation} & \bf $\bm{p(j)}$ & \bf support $\bm{S}$
 & $\bm{q(j)}$ & \bf parameters & \bf mean $\bm{\mu}$ & \XB{${\ds\textbf{classification}\atop\ds\textbf{rule}}$}\\
 \midrule
 \XA{{\bf1.} Poisson-type}{$X\sim P(\lambda)$} & $\ds e^{-\lambda}\frac{\lambda^j}{j!}$ & $\NN$ & $\lambda$ & $\lambda>0$ & $\lambda$ & $\delta=\beta=0$ \\
 [7ex]
 \XA{{\bf2.} Binomial-type}{$X\sim\mathrm{Bin}(N,p)$} & $\ds {N\choose{j}}p^j(1-p)^{N-j}$ & $0,1,\ldots,N$ & $p(N-j)$ & ${\ds N=1,2,\ldots\atop\ds 0<p<1}$
  & $Np$ & \XC{$\delta=0$\\$\beta=-p\in(-1,0)$} \\
 [7ex]
 \XA{{\bf3.} Negative Binomial-type}{$X\sim\mathrm{NB}(r,p)$} & $\ds \frac{[r]_j}{j!}p^r(1-p)^{j}$ & $\NN$ & $\ds \frac{1-p}{p}(r+j)$ & ${\ds r>0\atop\ds 0<p<1}$
  & $\ds \frac{r(1-p)}{p}$ & \XC{$\delta=0$\\$\ds \beta=\frac{1-p}{p}>0$} \\
 [7ex]
 \XA{\mbox{{\bf4a.} Negative} Hypergeometric-type}{$X\sim\mathrm{NHgeo}(N;r,s)$} & $\ds {N\choose{j}}\frac{(-r)_j(-s)_{N-j}}{(-r-s)_N}$ & $0,1,\ldots,N$
  & $\ds \frac{(r+j)(N-j)}{r+s}$ & ${\ds N=1,2,\ldots\atop \ds r,s>0}$ & $\ds \frac{Nr}{r+s}$ & $\ds \delta=\frac{-1}{r+s}<0$ \\
 \\
 [-1ex]
 \XA{\mbox{\bf4b.} Hypergeometric-type}{$X\sim\mathrm{Hgeo}(N;r,s)$} & $\ds {N\choose{j}}\frac{(r)_j(s)_{N-j}}{(-r-s)_N}$ & $0,1,\ldots,N$ & $\ds \frac{(r-j)(N-j)}{r+s}$
 & ${\ds N=1,2,\ldots\atop\ds r,s>N-1}$ & $\ds \frac{Nr}{r+s}$ & \XC{$\delta=\ds \frac{1}{r+s}>0$\\finite $S$} \\
 [7ex]
 \XA{{\bf5.} Discrete $F$-type\footnotemark\label{fnm:1}}{$X\sim\textrm{d-}F(\rho;r,s)$}
 & $\ds \frac{\varGamma(\rho-r)\varGamma(\rho-s)}{\varGamma(\rho-r-s)\varGamma(\rho)}
 \frac{[r]_j[s]_j}{j![\rho]_j}$
 & $\NN$ & $\ds \frac{(r+j)(s+j)}{\rho-r-s-1}$ & $
 {\ds (r,s)\in\CCC_2\atop\ds \rho>\max\{0,r+s+1\}
 \quad}$
 & $\ds \frac{rs}{\rho-r-s-1}$ & \XC{$\ds \delta=\frac{1}{\rho-r-s-1}>0$
 \\one-side infinite $S$} \\
 [7ex]
 \XA{{\bf6.} Discrete Student-type\footnotemark\label{fnm:2}$^,$\footnotemark\label{fnm:3}}{$X\sim\textrm{d-}t(\bm{z},\bm{w})$}
 & ${\ds C\frac{[z_1]_j[z_2]_j}{[w_1+1]_j[w_2+1]_j}}$
 & $\ZZ$
 & $\ds\frac{(z_1+j)(z_2+j)}{w_1+w_2-z_1-z_2}$ & $\ds{\bm{z},\bm{w}\in\tilde{\CCC}_2\colon\atop \ds w_1+w_2>z_1+z_2}$ & $\ds\frac{z_1z_2-w_1w_2}{w_1+w_2-z_1-z_2}$ & \XC{$\ds\delta=\frac{1}{w_1+w_2-z_1-z_2}>0$\\two-side infinite $S$} \\
 \bottomrule
 \end{tabular*}
 \bigskip

 \footnoterule
 \scriptsize
 \begin{tablenotes}
  \item[\ref{fnm:1}] $\CCC_2\doteq\{(z,\bar{z})\colon z\in\CC\smallsetminus\RR\}
                  \cup(0,\infty)^2\cup\left\{\bigcup_{n=0}^{\infty}(-n-1,-n)^2\right\}$.
  \item[\ref{fnm:2}] We were not able to find a closed formula for the normalizing constant $C$ of this pmf.
  \item[\ref{fnm:3}] $\tilde{\CCC}_2\doteq\{(z,\bar{z})\colon z\in\CC\smallsetminus\RR\}\cup\left\{\bigcup_{n\in\ZZ}(n,n+1)^2\right\}$.
  \end{tablenotes}
  \end{threeparttable}
 \end{sidewaystable}

 \section{A comparison with \citeauthor{Ord1}'s Discrete Student distributions}
 \label{section comparison}

 Here, we offer a comparison between
 \citeauthor{Ord1}'s Discrete Student distributions
 and the Discrete Student-type distributions
 that are presented in this article.

 \citet{Ord1} defined the Discrete Student distribution as one with pmf
 \begin{equation}
 \label{eq OrdDt}
 p(j)\propto \prod_{r=0}^k\frac{1}{(j+r+a)^2+b^2},\quad j\in\ZZ,
 \end{equation}
 where $k\in\NN$, $a\in[0,1]$ and $0<b^2<\infty$ are the parameters of the distribution.

 We are interested in answering the following questions:
 (a) Does $p$ in \eqref{eq OrdDt} belong to the CO system? (b) If yes, what is the corresponding pmf in the \Cref{table probabilities}?
 (c) Does \Cref{eq OrdDt} describe the set of Discrete Student-type distributions?

 Before our analysis, we state the following relations that arise
 by the definition of $[z]_j$; see \Cref{notations}\eqref{notations(c)}.
 Let $z\in\CC$, $r\in\NN^*$ and $j\in\ZZ$.
 Then, one can easily check that the following identities hold:
 \begin{equation}
 \label{eq [z]_j}
 [z]_{-j}=\left.{(-1)^j}\right/{[-z+1]_j}
 \quad\textrm{and}\quad
 [z]_{r+j}=[z]_j[z+j]_r=[z]_r[z+r]_j,
 \end{equation}
 provided that the quantities that appear are well-defined.

 Now, set the complex numbers $z_{1,2}=a\pm \imath b$, $w_{1,2}=a+k\pm \imath b$
 and then $\bm{z}_{k,a,b}=(z_1,z_2)$ and $\bm{w}_{k,a,b}=(w_1,w_2)$.
 Obviously, $\bm{z}_{k,a,b},\bm{w}_{k,a,b}\in\tilde{\CCC}_2$
 and $w_1+w_2-z_1-z_2=2k\in2\NN$.
 We observe that $\prod_{r=0}^k\left[(j+r+a)^2+b^2\right]=[z_1+j]_{k+1}[z_2+j]_{k+1}$.
 An application of \eqref{eq [z]_j} shows that $\prod_{r=1}^k\left[(j+r+a)^2+b^2\right]=[z_1]_{k+1}[z_2]_{k+1}[w_1+1]_j[w_2+1]_j/([z_1]_j[z_2]_j)$.
 Since the quantity $[z_1]_{k+1}[z_2]_{k+1}$ is positive and independent of $j$,
 the pmf in \eqref{eq OrdDt} takes the form
 \[
 p(j)\propto\frac{[z_1]_j[z_2]_j}{[w_1+1]_j[w_2+1]_j},\quad j\in\ZZ.
 \]
 If $k=0$, the pmf $p$ does not belong to the CO system;
 this is an expected result because $p$ does not have expected value
 due to the divergence of the harmonic series.
 In this case, $p$ is the pmf of a discrete Cauchy distribution.
 If $k\in\NN^*$, the pmf $p$ belongs to the CO system.
 In conclusion, let us denote the distribution of $p$
 in \eqref{eq OrdDt} by $\mathrm{d}^\textsc{ord}\textrm{-}t(k,a,b)$;
 then, for each $k\in\NN$ and $a,b$ as above,
 the pmf $p\sim\mathrm{d}^\textsc{ord}\textrm{-}t(k,a,b)$ belongs to the CO system
 iff $k\in\NN^*$; in particular,
 $\mathrm{d}^\textsc{ord}\textrm{-}t(k,a,b)\equiv\textrm{d-}t(\bm{z}_{k,a,b},\bm{w}_{k,a,b})$.

 Let $k\in\NN^*$, $0\le a\le1$ and $b^2>0$,
 and let us consider $p\sim \mathrm{d}^\textsc{ord}\textrm{-}t(k,a,b)$.
 Then, from the above analysis, it follows that
 $p\sim\textrm{d-}t(\bm{z}_{k,a,b},\bm{w}_{k,a,b})$.
 Observe that $\bm{z}_{k,a,b},\bm{w}_{k,a,b}\in\{(z,\bar{z})\colon z\in\CC\smallsetminus\RR\}\subsetneqq\tilde{\CCC}_2$
 with $w_1+w_2-z_1-z_2\in2\NN^*$.
 In view of \Cref{table probabilities}, it is obvious that
 the class of the Discrete Student-type distributions
 of the CO system is strictly bigger than \citeauthor{Ord1}'s class of the
 Discrete Student distributions which have finite mean value.
 Consequently, \eqref{eq OrdDt} cannot describe the whole of
 Discrete Student-type distributions of the CO system.
 Of course, we must note that the \citeauthor{Ord1}'s class of
 Discrete Student distributions contains discrete Cauchy distributions (for $k=0$);
 in contrast to \citeauthor{Ord1}'s system,
 any pmf of the CO system has finite mean value.

 Finally, it is worth noting the relationship between the finite moment-order
 and the parameter $k$ of $p\sim\mathrm{d}^\textsc{ord}\textrm{-}t(k,a,b)$.
 Consider the case $k=0$. Then, $p(j)\propto \left[(j+a)^2+b^2\right]^{-1}$
 and it is obvious that $p$ has finite moment of order $\theta$ iff $0\le\theta<1$.
 Suppose now that $k\in\NN^*$.
 Based on the previous analysis and \Cref{table probabilities},
 $p\sim\CO$ with $\delta=1/(2k)>0$.
 Since $1+1/\delta=2k+1$, \Cref{lem moments and delta} shows that
 $p$ has finite moment of order $\theta$ iff $0\le\theta<2k+1$.
 Observe that the rule
 ``$p$ has finite moment of order $\theta$ iff $0\le\theta<2k+1$''
 holds for every $k\in\NN$.

 \section{The symmetric pmfs of the CO system}
 \label{section symmetric}

 In this section, we are interested in characterizing the symmetric pmfs of the CO system.
 In investigating this aspect, we state and prove the following theorem.
 First, observe that if $X$ is a symmetric integer-valued rv with finite mean value,
 then the expected value of $X$ is an integer or half-integer number
 (the set of the half-integer numbers is denoted by $\ZZ+1/2$).

 \begin{theorem}
 \label{theo symmetric}
 Let $p\sim \CO$. The pmf $p$ is symmetric, around its mean value $\mu$,
 iff $\mu\in\frac{1}{2}\ZZ\doteq\ZZ\cup\{\ZZ+1/2\}$ and $4\delta\mu+2\beta=-1$.
 \end{theorem}

 \begin{proof}
 Suppose $X\sim p$. We prove separately the cases $\mu\in\ZZ$ and $\mu\in\ZZ+1/2$.

 Let $\mu\in\ZZ$ and let us consider the rv $Y=X-\mu$.
 Then, the rv $X$, and so the pmf $p$, is symmetric iff $Y\xlongequal{\mathrm{d}}-Y$.
 Using  \Cref{lem S:q:S_o:S^o}\eqref{lem S:q:S_o:S^o(d)},\eqref{lem S:q:S_o:S^o(e)},
 it follows that $Y\sim\mathrm{CO}(\mu_{Y};q_{Y})$, where $\mu_{Y}=0$
 and $q_{Y}(j)=\delta j^2+(2\delta\mu+\beta)j+\left(\delta\mu^2+\beta\mu+\gamma\right)$,
 and $-Y\sim\mathrm{CO}(\mu_{-Y};q_{-Y})$, where $\mu_{-Y}=0$
 and $q_{-Y}(j)=\delta j^2-(2\delta\mu+\beta+1)j+\left(\delta\mu^2+\beta\mu+\gamma\right)$.
 Applying \Cref{prop the support is defined by mu and q}\eqref{prop the support is defined by mu and q(b)},
 $Y\xlongequal{\mathrm{d}}-Y$ iff $4\delta\mu+2\beta=-1$.

 Let $\mu\in\ZZ+1/2$, say $\mu=\lfloor\mu\rfloor+1/2$.
 Consider the rvs $Y_1=X-\lfloor\mu\rfloor$ and $Y_2=-Y_1+1$.
 Then, the rv $X$, and so the pmf $p$, is symmetric iff $Y_1\xlongequal{\mathrm{d}}Y_2$.
 Again, from \Cref{lem S:q:S_o:S^o}\eqref{lem S:q:S_o:S^o(d)},\eqref{lem S:q:S_o:S^o(e)},
 we get $Y_1\sim\mathrm{CO}(\mu_{1};q_{1})$, where $\mu_{1}=1/2$
 and $q_{1}(j)=\delta j^2+(2\delta\lfloor\mu\rfloor+\beta)j+\left(\delta\lfloor\mu\rfloor^2+\beta\lfloor\mu\rfloor+\gamma\right)$,
 and $Y_2\sim\mathrm{CO}(\mu_{2};q_{2})$, where $\mu_{2}=1/2$
 and $q_{2}(j)=\delta j^2-[2\delta(\lfloor\mu\rfloor+1)+\beta+1]j+\left[\delta\lfloor\mu\rfloor^2+\delta(2\lfloor\mu\rfloor+1)+\beta(\lfloor\mu\rfloor+1)+\gamma+1/2\right]$.
 An application of \Cref{prop the support is defined by mu and q}\eqref{prop the support is defined by mu and q(b)} implies that
 $Y_1\xlongequal{\mathrm{d}}Y_2$ iff $4\delta\mu+2\beta=-1$.
 \end{proof}

 Now, we are interested in finding the types of the CO system that contain symmetric pmfs.
 If $X\sim \CO$, there exist $s\in\{-1,1\}$ and $r\in\ZZ$ such that
 the pmf of $Y=sX+r$ belongs to \Cref{table probabilities}.
 It is obvious that $X$ is a symmetric rv iff $Y$ is symmetric.
 Under this observation and using \Cref{theo symmetric} and \Cref{table probabilities},
 we have the following list:

  \begin{description}[itemsep=.5ex, wide, labelwidth=!, labelindent=0pt]
   \item[$\bullet$]
   The \textit{Poisson-type distributions} do not contain symmetric pmfs, due to non-symmetric support.
   Alternatively, since $\delta=\beta=0$, we have that $4\delta\mu+2\beta=0\ne-1$;

   \item[$\bullet$]
   The \textit{Binomial-type distributions} contain symmetric pmfs.
   Since $\delta=0$ and $\beta=-p$, $4\delta\mu+2\beta=-1$
   is equivalent with $p=1/2$ which implies $\mu=N/2\in\frac{1}{2}\ZZ$;

   \item[$\bullet$]
   The \textit{Negative Binomial-type} of distributions does not contain symmetric pmfs,
   due to non-symmetric support.
   Alternatively, $\delta=0$, $\beta=(1-p)/p$ and so $4\delta\mu+2\beta=2(1-p)/p>0$;

   \item[$\bullet$]
   The \textit{Negative Hypergeometric-type distributions} contain symmetric pmfs.
   If $p$ is a symmetric pmf, its support $S=\{0,1,\ldots,N\}$ must be symmetric around $\mu$;
   consequently, $\mu=rN/(r+s)=N/2$, equivalently, $r=s$.
   Conversely, for $r=s$, we have that $\mu=N/2$, $\delta=-1/(2r)$, $\beta=(N-r)/(2r)$ and so $4\delta\mu+2\beta=-1$;

   \item[$\bullet$]
   The \textit{Hypergeometric-type distributions} contain symmetric pmfs.
   Using the same arguments as in the Negative Hypergeometric-type distributions, a pmf in this subsystem is symmetric iff $r=s$;

   \item[$\bullet$]
   The \textit{Discrete $F$-type distributions} do not contain symmetric pmfs,
   due to non-symmetric support.
   Alternatively, setting $\theta=\rho-r-s-1>0$, we have that
   $\mu=rs/\theta$, $\delta=1/\theta$ and $\beta=(r+s)/\theta$. Observe that $rs>0$ and $r+s\in\RR$ because $(r,s)\in\CCC_2$.
   The relation $4\delta\mu+2\beta=-1$ implies that $\theta^2+2(r+s)\theta+4rs=0$, which has discriminant
   $\varDelta=4(s-r)^2$. If $r$ and $s$ are complex conjugate numbers,
   $\theta\in\CC\smallsetminus\RR$, a contradiction;
   therefore,
   $(r,s)\in(0,\infty)^2\cup\left\{\bigcup_{n=0}^{\infty}(-n-1,-n)^2\right\}$.+
   Solving the equation $\theta^2+2(r+s)\theta+4rs=0$, we get $\theta=-2r$ or $-2s$.
   Observe that $\mu=-s/2>0$ (or $\mu=-r/2>0$) belongs in $\frac{1}{2}\ZZ$.
   Hence, $s$ (or $r$) is a negative integer, a contradiction;

   \item[$\bullet$]
   The \textit{Discrete Student-type distributions} contain symmetric pmfs.
   Consider the vectors $\bm{z}=(-1/2,-1/2)$, $\bm{w}=(1/2,1/2)\in\tilde{\CCC}_2$.
   Then, $w_1+w_2-z_1-z_2=2$, $\mu=0$, $\beta=-1/2$ and so $4\delta\mu+2\beta=-1$.
   Using \eqref{eq pmf delta>0: S=Z}, we find that
   the corresponding symmetric pmf is
   $p(0)=(2\pi^2-55/3)^{-1}$, $p(j)=(6\pi^2-165)^{-1}$ when $j=\pm1$,
   and $p(j)=(2\pi^2-55/3)^{-1}(j^2-1/4)^{-2}$ for $j=\pm2,\pm3,\ldots$~.

   Now, we determine the class of the symmetric discrete-$t$ rvs.
   Suppose $X\sim\textrm{d-}t(\bm{z},\bm{w})$,
   where $\bm{z},\bm{w}\in\tilde{\CCC}_2$ with $\delta^{-1}=w_1+w_2-z_1-z_2>0$.
   Then, $X\sim\mathrm{CO}(\mu_X;q_X)$ for an admissible pair $(\mu_X;q_X)$,
   and the rv $Y=X-\lfloor\mu_X\rfloor$ follows $\mathrm{CO}(\mu_Y;q_Y)$,
   see \Cref{lem S:q:S_o:S^o}\eqref{lem S:q:S_o:S^o(d)},
   and has mean value $\mu_Y\in[0,1)$.
   Since $X$ is a symmetric rv iff $Y$ is symmetric, it is sufficient to find the generator class
   of the symmetric discrete-$t$ distributions for which the mean value is $0$ or $1/2$.
   We distinguish the cases $\mu=0$ and $\mu=1/2$.
  \begin{description}[leftmargin=10pt, labelwidth=!]
   \item[--]
   Case $\mu=0$. In view of \Cref{table probabilities}, $\mu_X=\delta(z_1z_2-w_1w_2)=0$, $q_X(j)=\delta(z_1+j)(z_2+j)$ and
   $p_X(j)\propto[z_1]_j[z_2]_j/([w_1+1]_j[w_2+1]_j)$, $j\in\ZZ$.
   The rv $Y=-X$ has pmf $p_Y(j)=p_X(-j)\propto[z_1]_{-j}[z_2]_{-j}/([w_1+1]_{-j}[w_2+1]_{-j})$, $j\in\ZZ$.
   Applying \eqref{eq [z]_j}, $p_Y(j)\propto[-w_1]_{j}[-w_2]_{j}/([-z_1+1]_{j}[-z_2+1]_{j})$, $j\in\ZZ$.
   \Cref{lem S:q:S_o:S^o}\eqref{lem S:q:S_o:S^o(d)} gives that $Y$ follows $\mathrm{CO}(\mu_Y;q_Y)$;
   moreover, from  \Cref{table probabilities}, $\mu_Y=\delta(w_1w_2-z_1z_2)=0$
   and $q_Y(j)=\delta(-w_1+j)(-w_2+j)$. Obviously, $X$ is a symmetric rv iff $X\xlongequal{\mathrm{d}}Y$;
   using \Cref{prop the support is defined by mu and q}\eqref{prop the support is defined by mu and q(b)},
   $X$ is a symmetric rv iff $\bm{w}\eqcirc-\bm{z}$ and $z_1+z_2<0$ (since $\delta>0$). The symmetric discrete-$t$ rvs
   with mean value zero is the set
   \[
   \CS_{\textrm{d-}t}^{\{0\}}\doteq\left\{\textrm{d-}t(\bm{z},-\bm{z})\colon \bm{z}\in\tilde{\CCC}_2\textrm{ with } z_1+z_2<0\right\};
   \]

   \item[--]
   Case $\mu=1/2$. Again from \Cref{table probabilities} we have that $\mu_X=\delta(z_1z_2-w_1w_2)=1/2$ and $q_X(j)$, $p_X(j)$ as in the case $\mu=0$.
   The rv $Y=-X+1$ has pmf $p_Y(j)=p_X(-j+1)\propto[z_1]_{-j+1}[z_2]_{-j+1}/([w_1+1]_{-j+1}[w_2+1]_{-j+1})$, $j\in\ZZ$.
   Applying \eqref{eq [z]_j},
   $p_Y(j)\propto[-w_1-1]_{j}[-w_2-1]_{j}([-z_1]_{j}[-z_2]_{j})$, $j\in\ZZ$.
   \Cref{lem S:q:S_o:S^o}\eqref{lem S:q:S_o:S^o(d)},\eqref{lem S:q:S_o:S^o(e)} give that $Y$ follows $\mathrm{CO}(\mu_Y;q_Y)$.
   By construction, $\mu_Y=\mu_X=1/2$; furthermore, \Cref{table probabilities} implies
   $q_Y(j)=\delta(-w_1-1+j)(-w_2-1+j)$. Since $X$ is a symmetric rv iff $X\xlongequal{\mathrm{d}}Y$,
   \Cref{prop the support is defined by mu and q}\eqref{prop the support is defined by mu and q(b)} gives that
   $X$ is a symmetric rv iff $\bm{w}\eqcirc-\bm{z}-\bm{1}$ and $z_1+z_2<-1$ (since $\delta>0$). The symmetric discrete-$t$ rvs
   with mean value half is
   \[
   \CS_{\textrm{d-}t}^{\{1/2\}}\doteq\left\{\textrm{d-}t(\bm{z},-\bm{z}-\bm{1})\colon \bm{z}\in\tilde{\CCC}_2\textrm{ with } z_1+z_2<-1\right\}.
   \]
  \end{description}
  The symmetric discrete-$t$ rvs with mean value zero or half is the generator-set of the symmetric discrete-$t$ rvs,
  \[
  \CS_{\textrm{d-}t}^{\{0,1/2\}}=\CS_{\textrm{d-}t}^{\{0\}}\cup\CS_{\textrm{d-}t}^{\{1/2\}}.
  \]
  The set of the symmetric discrete-$t$ rvs is
  \[
  \CS_{\textrm{d-}t}\doteq\CS_{\textrm{d-}t}^{\{0,1/2\}}+\ZZ=\left\{X+r\colon X\in\CS_{\textrm{d-}t}^{\{0,1/2\}}, \ r\in\ZZ\right\}.
  \]
 \end{description}

 Finally, we define the  noncentrality parameter as well as the degrees of freedom of a discrete-$t$ distribution.
 In view of \Cref{theo symmetric} and the fact that the $t_\nu$ distribution has finite absolute moments of order $\theta$
 for each $0<\theta<\nu$ while its $\nu$th absolute moment is infinity, cf.\ \Cref{lem moments and delta}, we give the following definition.

 \begin{definition}
 \label{def D.t}
 Let $X\sim\textrm{d-}t(\bm{z},\bm{w})$ with $\bm{z},\bm{w}\in\tilde{\CCC}_2$ and $\delta^{-1}=w_1+w_2-z_1-z_2>0$,
 and let us consider the parameters $\mu=\delta(z_1z_2-w_1w_2)$ and $\beta=\delta(z_1+z_2)$.
 Then, the \textit{noncentrality parameter} of $X$ is defined by
 \[
 \bm{\varepsilon}=(\varepsilon_1,\varepsilon_2),
 \quad\textrm{where} \
 \varepsilon_1=\min_{x\in\frac{1}{2}\ZZ}\left|x-\mu\right|
 \textrm{ and }
 \varepsilon_2=|4\delta\mu+2\beta+1|,
 \]
 and the \textit{degrees of freedom} of $X$ are defined as $\mathrm{df}=1+1/\delta$.
 \end{definition}

 \section{Moment relations in the Cumulative Ord family}
 \label{section moment relations in CO}
 This section presents some properties about the moments of a
 rv of the CO family.

 For a discrete rv $X\sim\COq\equiv\CO$,
 the following covariance identity holds
 \begin{equation}
 \label{eq cov-ident}
 \Cov[X,g(X)]=\E[q(X)\Delta g(X)],
 \end{equation}
 provided that $\E[q(X)|\Delta g(X)|]<\infty$;
 see \cite{CP2}. Setting
 $g(x)=x$, we get
 \[
 \sigma^2\doteq\Var(X)=\E[q(X)],
 \]
 provided $\E\left(X^2\right)<\infty$. Writing
 $q(X)=\delta(X-\mu)^2+q'(\mu)(X-\mu)+q(\mu)$ and
 taking expectations, we have
 \begin{equation}
 \label{eq Var(X)=q(mu)/(1-delta)}
 \sigma^2={q(\mu)}/{(1-\delta)},
 \end{equation}
 noting that from \Cref{lem moments and delta} and
 \Cref{rem finite S: delta>0,rem |S|=2},
 the denominator $1-\delta$ is positive.

 Now, we prove a lemma concerning the pmf
 $p^*(j)\propto q(j)p(j)$.
 \begin{lemma}
 \label{lem q(j)p(j)}
 Suppose a non-constant rv $X\sim\COq$ and
 $\E|X|^3<\infty$. Let $X^*$ be the rv with pmf
 $p^*\propto qp$. Then, $X^*$ is supported
 on the set ${S(X^*)}\doteq S^\circ(X)$ (i.e., $\alpha^*=\alpha$,
 $\omega^*=\omega-1$) and $X^*\sim\mathrm{CO}(\mu^*;q^*)$, where
 $\mu^*={(\mu+\beta+\delta)}/{(1-2\delta)}$ and
 $q^*(j)={q(j+1)}/{(1-2\delta)}$.
 \end{lemma}
 \begin{proof}
 \Cref{lem moments and delta} proves that $1-2\delta>0$ because $\E|X|^3<\infty$;
 from this and \Cref{lem S:q:S_o:S^o}\eqref{lem S:q:S_o:S^o(b)},\eqref{lem S:q:S_o:S^o(c)},
 it follows that the function $p^*$ is non-negative on $\ZZ$ and takes strictly
 positive values on the set ${S(X^*)}=S^\circ(X)$.
 Using
 $\Delta[h_1(k)h_2(k)]=h_1(k)\Delta{h_2(k)}+h_2(k+1)\Delta{h_1(k)}$, we have
 $\Delta[q(k)q(k-1)p(k-1)]=q(k)\Delta[q(k-1)p(k-1)]+q(k)p(k)\Delta{q(k)}
                         =q(k)(\mu-k)p(k)+q(k)p(k)(2\delta{k}+\delta+\beta)=[\mu+\delta+\beta-(1-2\delta)k]q(k)p(k)$.
 Thus,
 $\sum_{k\le{j}}(\mu^*-k)p^*(k) =\left.\left\{\sum_{k\le{j}}\Delta[q(k)q(k-1)p(k-1)]\right\}\right/\{(1-2\delta)\E[q(X)]\}$.
 Since $\E|X|^3<\infty$, using the same arguments as in the
 proof of \Cref{prop CO against Ord system},
 we obtain that
 $\sum_{k\le{j}}\Delta[q(k)q(k-1)p(k-1)]=q(j+1)q(j)p(j)$.
 So, $\sum_{k\le{j}}(\mu^*-k)p^*(k)=q^*(j)p^*(j)$. It
 remains to show that the value $\mu^*$ is the mean of
 $X^*$. Of course, $\E|X^*|<\infty$ because $\E|X|^3<\infty$.
 If $S(X^*)$ has a finite upper endpoint
 $\omega^*=\omega-1$, then
 $q^*(\omega^*)={q(\omega)}/{(1-2\delta)}=0$, since
 $\omega$ is the upper endpoint of $S(X)$. If
 $\omega^*=\omega=\infty$, we use the same arguments as in the
 proof of \Cref{prop CO against Ord system}. For
 both cases $\omega<\infty$ and $\omega=\infty$,
 $\sum_{k\in{S}(X^*)}(\mu^*-k)p^*(k)=0$, i.e., $\E\left(X^*\right)=\mu^*$.
 \end{proof}

 The quadratic $q$ takes non-negative values on the support
 of $X$. Therefore, we can create new pmfs by defining
 $p_i\propto q^{[i]}p$. But, if the support of $X$ is finite,
 then for each $i$ greater than or equal to the cardinality
 of $S(X)$, the function $p_i$ vanishes identically on
 $\ZZ$. Thus, it is useful to define the quantity $M=M(X)$
 as follows:
 \[
 M=M(X)\doteq|S(X)|-1=\omega-\alpha\in\{0,1,\ldots\}\cup\{\infty\}.
 \]
 \begin{proposition}
 \label{prop q^[i](j)p(j)}
 Let $X\sim\COq=\CO$ with pmf $p$ and $\E|X|^{2n+1}<\infty$
 for some $n\in\{0,1,\ldots,M(X)\}$. For all $i=0,1,\ldots,n$, we consider the rvs
 $X_i$ with pmfs
 $p_i\propto q^{[i]}p$
 [note that $X_i=X_{i-1}^*$, $i=1,\ldots,n$, where $X_0=X$],
 and we define
 \begin{align*}
 &\mu_i=\frac{\delta{i^2}+\beta{i}+\mu}{1-2i\delta},
 &&q_i(j)=\frac{q(j+i)}{1-2i\delta},\\
 &\delta_i=\frac{\delta}{1-2i\delta},
 &&\psi(i)=\frac{q\left.\left(\left[-\delta{i^2}+(\beta+1)i+\mu\right]\right/[1-2i\delta]\right)}{1-(2i+1)\delta}.
 \end{align*}
 Then:
 \begin{enumerate}[itemsep=.5ex, wide, labelwidth=!, labelindent=0pt, label=\rm(\alph*), ref=\textcolor{black}{\alph*}]
 \item
 \label{prop q^[i](j)p(j)(a)}
 The rv $X_i$ is supported on the set
 $S_i\doteq S(X_i)=\{\alpha_i,\ldots,\omega_i\}=\{\alpha,\ldots,\omega-i\}$;
 \item
 \label{prop q^[i](j)p(j)(b)}
 $X_i\sim\mathrm{CO}(\mu_i;q_i)$;
 \item
 \label{prop q^[i](j)p(j)(c)}
 $\Var(X_i)=\psi(i)$
 (for $i=n$, it is additionally required that
 $\E|X|^{2n+2}<\infty$);
 \item
 \label{prop q^[i](j)p(j)(d)}
 $A_i=A_i(\mu;q)\doteq\E\left[q^{[i]}(X)\right]=\psi^{[i]}(0)\Pi_{2\delta}^{[i]}(0)$;
 \item
 \label{prop q^[i](j)p(j)(e)}
 The descending factorial moments of $X$, $\mu_{(r)}=\E(X)_r$,
 and the ascending factorial moments of $X$, $\mu_{[r]}=\E[X]_r$,
 satisfy the following second-order recurrence relations:
 \[
 (1-r\delta)\mu_{(r+1)}=\{\mu+r[\beta-1+(2r-1)\delta]\}\mu_{(r)}+r\{\gamma+(r-1)[\beta+(r-1)\delta]\}\mu_{(r-1)},
 \]
 \[
 (1-r\delta)\mu_{[r+1]}=\{\mu+r[\beta+2-(2r-1)\delta]\}\mu_{[r]}+r\{\gamma-\mu-(r-1)[\beta+1-(r-1)\delta]\}\mu_{[r-1]},
 \]
 with initial conditions $\mu_{(0)}=\mu_{[0]}=1$ and $\mu_{(1)}=\mu_{[1]}=\mu$, for all $r=1,\ldots,2n$;
 \item
 \label{prop q^[i](j)p(j)(f)}
 The factorial moments of $X$,
 $\mu_{(r)}$ and $\mu_{[r]}$,
 satisfy the following recurrence relations:
 \[
 (1-r\delta)\mu_{[r+1]}=[\mu+r(\beta-r\delta+1)]\mu_{[r]}+\gamma r!\sum_{k=0}^{r-1}\frac{\mu_{[k]}}{k!},
 \]
 \[
 (1-r\delta)\mu_{(r+1)}=\left(\delta r^2+\beta r+\mu\right)\mu_{(r)}+(\gamma-\mu) r!\sum_{k=0}^{r-1}(-1)^{r+k+1}\frac{\mu_{(k)}}{k!}.
 \]
 \end{enumerate}
 \end{proposition}
 \begin{proof}
\eqref{prop q^[i](j)p(j)(a)}
Observe that $1-2i\delta>0$ for all $i=0,1,\ldots,n$
(if $\delta\le0$, it is obvious; if $\delta>0$,
the case of infinity support follows by \Cref{lem moments and delta}
while the case of finite support by \Cref{rem finite S: delta>0}).
Therefore, \Cref{lem S:q:S_o:S^o}\eqref{lem S:q:S_o:S^o(b)},\eqref{lem S:q:S_o:S^o(c)}
show that $q^{[i]}p$ is supported on $S_i$.

\eqref{prop q^[i](j)p(j)(b)}
The proof will be done by induction on $i$. For $i=1$,
the result follows from \Cref{lem q(j)p(j)}. Assuming
that it holds for $i-1\in\{0,1,\ldots,n-1\}$, we will prove
that it is true for $i$. By assumption, $X_{i-1}\sim\mathrm{CO}\left(\mu_{i-1};q_{i-1}\right)$.
From $\E|X|^{2n+1}<\infty$, it follows that
$\E|X_{i-1}|^{3}<\infty$. As in \Cref{lem q(j)p(j)}, we
consider the rv $X^*_{i-1}\sim\mathrm{CO}\left(\mu^*_{i-1};q^*_{i-1}\right)$ with
$\mu^*_{i-1}=(\mu_{i-1}+\beta_{i-1}+\delta_{i-1})/(1-2\delta_{i-1})$
and $q^*_{i-1}=q_{i-1}(j+1)/(1-2\delta_{i-1})$.
Hence, after some algebra, we get $\mu^*_{i-1}=\mu_i$ and
$q^*_{i-1}=q_i$. Finally, observe that
$p^*_{i-1}\propto q_{i-1}p_{i-1}$ and $q_{i-1}p_{i-1}\propto q^{[i]}p$;
so, $p^*_{i-1}\propto q^{[i]}p$.
By definition, $p_i\propto q^{[i]}p$.
Hence, we conclude that $p_i=p^*_{i-1}$ because $p_i$, $p^*_{i-1}$ are pmfs with support $S_i$.

\eqref{prop q^[i](j)p(j)(c)}
It is immediate from \eqref{prop q^[i](j)p(j)(b)} and \eqref{eq Var(X)=q(mu)/(1-delta)}.

\eqref{prop q^[i](j)p(j)(d)}
 Using \eqref{prop q^[i](j)p(j)(c)}, an application of \eqref{eq Var(X)=q(mu)/(1-delta)} gives
 $(1-2j\delta)\psi(j)=(1-2j\delta)\E\left[q_j(X_j)\right]=\E\left[q(X_j+j)\right]={\E\left[q^{[j+1]}(X)\right]}\left/{\E\left[q^{[j]}(X)\right]}\right.={A_{j+1}}/{A_{j}}$, $j=0,1,\ldots,n-1$,
 where $A_0=1$. By multiplying these
 relations for $j=0,1,\ldots,k-1$, the result follows.

 \eqref{prop q^[i](j)p(j)(e)}
 Write
 $\E(X)_{r+1}=(\mu-r)\E(X)_{r}+\E[(X-\mu)(X)_{r}]=(\mu-r)\E(X)_{r}+\Cov[X,(X)_{r}]$.
 Using the covariance identity \eqref{eq cov-ident} and since
 $\Delta(j)_{r}=r(j)_{r-1}$, it follows that
 $\Cov[X,(X)_{r}]=r\E[q(X)(X)_{r-1}]$.
  Moreover,
  $q(x)(x)_{r-1}=\delta(x)_{r+1}+[\beta+(2r-1)\delta](x)_{r}
  +\{\gamma+(r-1)[\beta+(r-1)\delta]\}(x)_{r-1}$.
 Thus,
 $\E[q(X)(X)_{r-1}]=\delta\mu_{(r+1)}+[\beta+(2r-1)\delta]\mu_{(r)}
 +\{\gamma+(r-1)[\beta+(r-1)\delta]\}\mu_{(r-1)}$.
 Upon combining the above relations, the result follows.

 For the second relation, we consider the rv
 $Y=-X\sim\mathrm{CO}(-\mu;\delta,-\beta-1,\gamma-\mu)$;
 see \Cref{lem S:q:S_o:S^o}\eqref{lem S:q:S_o:S^o(e)}.
 Observe that $\mu_{(r)}^Y=(-1)^{r}\mu_{[r]}$.
 An application of the first relation, with some algebra, shows the result.

 \eqref{prop q^[i](j)p(j)(f)}
 Using the same arguments as in \eqref{prop q^[i](j)p(j)(e)}, we write
 $\mu_{[r+1]}=(\mu+r)\mu_{[r]}+\Cov(X,[X]_{r})$. Utilizing
 \eqref{eq cov-ident} and the fact that
 $\Delta[X]_{r}=r[X+1]_{r-1}$,
 we get $\Cov(X,[X]_{r})=r\E\{q(X)[X+1]_{r-1}\}$.
 Write
 $q(x)[x+1]_{r-1}=\delta[x]_{r+1}+(\beta-r\delta)[x]_{r}+\gamma\sum_{k=0}^{r-1}(r-1)_{r-1-k}[x]_k$.
 Then,
 $\mu_{[r+1]}=(\mu+r)\mu_{[r]}+r\left[\delta\mu_{[r+1]}+(\beta-r\delta)\mu_{[r]}
 +\gamma(r-1)!\sum_{k=0}^{r-1}\mu_{[k]}/k!\right]$.
 Finally, the proof of the recurrence relation of $\mu_{(r)}$s is similar to that of \eqref{prop q^[i](j)p(j)(e)}.
 \end{proof}

 Now, suppose the rv $X$ belongs to the CO family and its support has lower endpoint $\alpha=0$.
 Then, $\gamma=\mu$ (see \Cref{lem S:q:S_o:S^o}\eqref{lem S:q:S_o:S^o(c)}) and so
 the second recurrence relation
 of \Cref{prop q^[i](j)p(j)}\eqref{prop q^[i](j)p(j)(f)}
 takes the form $(1-r\delta)\mu_{(r+1)}=q(r)\mu_{(r)}$.
 Under this observation, the following corollary follows immediately.

 \begin{corollary}
 \label{cor recurrence}
 Let $X\sim \CO$. If the support of $X$ has lower endpoint $\alpha=0$,
 then for each positive integer $k$ such that $\E|X|^k<\infty$, the $k$th
 descending factorial moment of $X$ is
 $\mu_{(k)}={q^{[k]}(0)}\left/{\Pi^{[k]}_{\delta}(0)}\right.=\prod_{j=0}^{k-1}[q(j)/(1-j\delta)]$.
 \end{corollary}

 We apply \Cref{cor recurrence} to the distributions of the types 1--5 that are presented in \Cref{table probabilities}.

 \begin{application}
 \label{application}
 \begin{description}[itemsep=.5ex, wide, labelwidth=!, labelindent=0pt]
 \item[{\small 1.} \sc Poisson distribution:]
 If $X\sim P(\lambda)$ with $\lambda>0$, then $q(j)/(1-j\delta)=\lambda$
 and so $\mu_{(k)}=\lambda^k$ for all $k=0,1,\ldots$;
 \item[{\small 2.} \sc Binomial distribution:]
 If $X\sim\mathrm{Bin}(N,p)$ with $N=1,2,\ldots$ and $0<p<1$, then $q(j)/(1-j\delta)=p(N-j)$
 and so $\mu_{(k)}=p^k(N)_k$ for all $k=0,1,\ldots$;
 \item[{\small 3.} \sc Negative Binomial distribution:]
 If $X\sim\mathrm{NB}(r,p)$, $r>0$ and $0<p<1$, then $q(j)/(1-j\delta)=[(1-p)/p](r+j)$
 and so $\mu_{(k)}=[(1-p)/p]^k[r]_k$ for all $k=0,1,\ldots$;
 \item[{\small 4a.} \sc Negative Hypergeometric distribution:]
 If $X\sim\mathrm{NHgeo}(N;r,s)$ with $N=1,2,\ldots$ and $r,s>0$, then $q(j)/(1-j\delta)=(r+j)(N-j)/(r+s+j)$
 and so $\mu_{(k)}=[r]_k(N)_k/[r+s]_k$ for all $k=0,1,\ldots$;
 \item[{\small 4b.} \sc Hypergeometric distribution:]
 If $X\sim\mathrm{NHgeo}(N;r,s)$ with $N=1,2,\ldots$ and $r,s>N-1$, then $q(j)/(1-j\delta)=(r-j)(N-j)/(r+s-j)$
 and so $\mu_{(k)}=(r)_k(N)_k/(r+s)_k$ for all $k=0,1,\ldots$;
 \item[{\small 5.} \sc Discrete $F$-type distribution:]
 If $X\sim\textrm{d-}F(\rho;r,s)$ with $(r,s)\in\CCC_2$ and $\rho>\max\{0,r+s+1\}$, then $q(j)/(1-j\delta)=(r+j)(s+j)/(\rho-r-s-1-j)$
 and so $\mu_{(k)}=[r]_k[s]_k/(\rho-r-s-1)_k$ for all $k=0,1,\ldots$ such that $k<\rho-r-s$.
  \end{description}
 \end{application}

 Next, we generalize the results of \Cref{lem lim_j->oo q(j)p(j)} in CO family.

\begin{proposition}
\label{prop lim_j->oo j^2ip(j)}
Let $X\sim\COq=\CO$ and assume that it has an upper (resp.\
lower) unbounded support and $\E|X|^{2i-1}<\infty$ for some
$i\in\{1,2,\ldots\}$. Then, $j^{2i}p(j)\to0$ as $j\to
\infty$ (resp.\ $j\to-\infty$).
 \end{proposition}
 \begin{proof}
Note that $\delta\geq 0$ since
 the support is finite if $\delta<0$.
For the case $\delta=0$, the result is obvious since $X$
has finite moments of any order; see \Cref{lem moments
and delta}. When $\delta>0$, then, as in \Cref{prop q^[i](j)p(j)},
consider the rv $X_{2i-2}\sim\mathrm{CO}(\mu_{2i-2};q_{2i-2})$.
From \Cref{lem lim_j->oo q(j)p(j)}, it follows that $q_{2i-2}(j)p_{2i-2}(j)\to0$ as
$j\to \infty$ (resp.\ $j\to-\infty$) and since
$\lim_{j\to\pm\infty}q_{2i-2}(j)p_{2i-2}(j)\propto\lim_{j\to\pm\infty}j^{2i}p(j)$,
the proof is complete.
\end{proof}

 \section{Orthogonal polynomials in the Cumulative Ord family}
 \label{section orthogonal polynomials in CO}
 In this section, we
 present results for the orthogonal polynomials of a probability
 measure of the CO family. These polynomials are obtained
 by a {\it discrete Rodrigues-type formula}.

 First, we present a brief review.
 \citet[Chap.~IV, pp.~419--439]{Hild} studied the
 nonzero solutions $u(j)$ of the
 Pearson difference equation,
 \begin{equation}
 \label{eq Pearson Diff-eq}
 \Delta u(j)=\frac{N(j)}{D(j)}u(j),
 \end{equation}
 where the numerator $N$ is a polynomial of degree at most
 one and the denominator $D$ is a polynomial of degree at
 most two. He showed that the functions $Q_n(j)$, produced by the
 Rodrigues-type formula
 \begin{equation}
 \label{eq Rodiquez Hild}
 Q_n(j)=\frac{\Delta^n\left[D^{[n]}(j-n)u(j)\right]}{u(j)},
 \end{equation}
 are polynomials of degree at most $n$; see
 \cite[p.~425]{Hild}.
 Note that \citeauthor{Hild} makes use of the
 descending power notation,
 $D^{(n)}(j-1)\doteq D(j-1)D(j-2)\cdots{D(j-n)}=D^{[n]}(j-n)$.
 He farther established several properties of these polynomials.
 In the sequel of this section, when we say that a function
 is the solution of a difference equation, we will always
 mean a pmf solution.

 In \citeauthor{Hild}'s results, the orthogonality
 of the produced polynomials was not an issue.
 However, these polynomials are orthogonal
 only when we make a correct choice of the set on which we
 seek a solution, and provided that we used
 the correct writing of the ratio of the polynomials
 $N$ and $D$ in \eqref{eq Rodiquez Hild}.
 Next, we present some examples to illustrate this issue.

 Here, we note that the \cref{eq Ord's Diff-eq,eq Pearson Diff-eq} are equivalent, excluding
 the case $\Delta p(j)=0$. Specifically,
 ${\Delta p(j)}/{p(j)}={N(j)}/{D(j)}$
 is equivalent with
 ${\Delta p(j-1)}/{p(j)}={N(j-1)}/{[D(j-1)+N(j-1)]}$.

  \begin{example}
 \label{exm Poisson:Uniform:Geometric}
 \begin{enumerate}[itemsep=.5ex, wide, labelwidth=!, labelindent=0pt, label=\rm(\alph*), ref=\textcolor{black}{\alph*}]
 \item
 \label{exm Poisson:Uniform:Geometric(a)}
 Consider the difference equation
 ${\Delta{p(j)}}/{p(j)}={(\lambda-j-1)}/{(j+1)}$,
 where $\lambda$ is a positive constant. This difference
 equation is of the
 form \eqref{eq Ord's Diff-eq} and \eqref{eq Pearson Diff-eq}.
 Of course, in order to
 solve a difference equation, we must specify the support set on
 which we seek the solution. If this set is $\NN$, then the
 solution is $\left.e^{-\lambda}{\lambda^j}\right/{j!}$,
 $j=0,1,\ldots$ (Poisson distribution with parameter
 $\lambda$). If the set is $\{0,1,\ldots,N\}$, then the
 solution is $\left.C{\lambda^j}\right/{j!}$, $j=0,1,\ldots,N$
 (truncated Poisson distribution with parameter $\lambda$).
 The polynomials obtained by \eqref{eq Rodiquez Hild}
 are the  Charlier polynomials which are orthogonal
 with respect to the Poisson pmf, but not
 with respect to the truncated Poisson pmf.
 \item
 \label{exm Poisson:Uniform:Geometric(b)}
 Consider the pmf of the geometric distribution with parameter
 $p\in(0,1)$, i.e., $p(j)=p(1-p)^j$, $j=0,1,\ldots$~. This
 pmf  satisfies the difference equation ${\Delta
 p(j)}/{p(j)}=-p$, which can be rewritten in the form
 \eqref{eq Pearson Diff-eq} in many ways. Specifically,
 ${\Delta p(j)}/{p(j)}={-p(bj+a)}/{(bj+a)}$,
 where $bj+a$ is a constant (when $b=0$), or a linear polynomial
 without roots on $\NN$. For any choice of $a$
 and $b$, \citeauthor{Hild}'s results are valid. However, the
 polynomials in \eqref{eq Rodiquez Hild} are orthogonal with
 respect to the geometric pmf only when we make the choice $a=b\ne0$ (Meixner
 polynomials).
 \item
 \label{exm Poisson:Uniform:Geometric(c)}
 Now, consider the difference equation
 ${\Delta{p(j)}}/{p(j)}=0$ supported on an
 integer chain. Of course, if the support is infinite, then
 it has no pmf solutions; thus, we consider a finite integer
 chain, and without loss of generality take $S=\{1,2,\ldots,N\}$.
 The solution is $p(j)={1}/{N}$, $j=1,2,\ldots,N$, i.e.,
 $X$ is uniformly distributed on the support. The equation can be
 rewritten in the form \eqref{eq Pearson Diff-eq} in
 many ways, i.e., $N(j)=0$ and $D(j)=cj^2+bj+a$ any quadratic
 polynomial without roots on $\{1,2,\ldots,N-1\}$. Again, the
 polynomials in \eqref{eq Rodiquez Hild} are orthogonal with
 respect to pmf $p$ only when one makes the correct choice
 $D(j)\propto
 j(N-j)$ (Hahn polynomials).
 \end{enumerate}
 \end{example}

 It is true that the denominator in \eqref{eq Ord's Diff-eq},
 under suitable conditions, generates orthogonal polynomials
 with respect to the pmf solution of this equation; see
 \Cref{prop CO against Ord system} and also
 the next theorem.

 \begin{remark}
  \label{rem for Example}
 In view of \Cref{exm Poisson:Uniform:Geometric}, we observe the following.
 The Rodrigues-type formula \eqref{eq Rodiquez Hild} is a mechanism for  producing  polynomials,
 that may have some elegant properties regarding their coefficients.
 On the other hand, the specific cases of \Cref{exm Poisson:Uniform:Geometric} clearly indicate that the relation \eqref{eq Ord's Diff-eq}
 (or the equivalent relation \eqref{eq Pearson Diff-eq}) ignores the information about the production of the
 Rodrigues-orthogonal polynomials, while the relation \eqref{eq CO}
 provides the whole of the information that is needed.
 \end{remark}

 Independently of \citeauthor{Hild}'s results, \citet{APP2} studied
 the orthogonality of the Rodrigues polynomials in the CO
 family:
 \begin{theorem}[\textrm{\cite[Lemma 2.3, Theorems 2.1 and 2.2]{APP2}}]
 \label{theo orthogonal polynomials in CO}
 Let $X\sim\COq=\CO$. For each $k=0,1,2,\ldots$,
 define the functions $P_k(j)$, $j\in S$, by the
 Rodrigues-type formula
 \begin{equation}
 \label{eq Rodiquez APP}
 P_k(j)=\frac{(-1)^k}{p(j)}\Delta^k\left[q^{[k]}(j-k)p(j-k)\right]
       =\frac{1}{p(j)}\sum_{i=0}^{k}(-1)^{k-i}{k\choose{i}}q^{[k]}(j-i)p(j-i).
 \end{equation}
 Then:
 \begin{enumerate}[itemsep=.5ex, wide, labelwidth=!, labelindent=0pt, label=\rm(\alph*), ref=\textcolor{black}{\alph*}]
 \item
 \label{theo orthogonal polynomials in CO(a)}
 Each $P_k$ is a polynomial of degree at most $k$, with
 \begin{equation}
 \label{eq lead P_k}
 \lead(P_k)=\Pi_\delta^{[k]}(k-1)\doteq c_k(\delta)
 \end{equation}
 {\rm[in the sense that the function $P_k(j)$, $j\in S$, is
 the restriction of a real polynomial
 $G_k(x)=\sum_{i=0}^{k}c(k,i)x^i$, $x\in\RR$, of degree at
 most $k$, such that $c(k,k)=\lead(P_k)$]};
 \item
 \label{theo orthogonal polynomials in CO(b)}
 Provided that $\E|X|^{2n}<\infty$ for some $n\ge1$,
 the polynomials $P_k$, $k=0,1,\ldots,n$, satisfy the
 orthogonality condition
 \begin{equation}
 \label{eq orthogonality condition}
 \E[P_k(X)P_m(X)]=\delta_{k,m}c_k(\delta)
 \E\left[q^{[k]}(X)\right]=\delta_{k,m}c_k(\delta)A_k,
 \quad
 k,m=0,1,\ldots,n,
 \end{equation}
 where $\delta_{k,m}$ is Kronecker's delta;
 \item
 \label{theo orthogonal polynomials in CO(c)}
 Provided that $k\ge1$ and $\E|X|^{2k-1}<\infty$,
 the following ``Rodrigues inversion formula'' holds:
 \begin{equation}
 \label{eq Rodrigues inversion formula}
 q^{[k]}(j)p(j)=\frac{1}{(k-1)!}\sum_{i>j}(i-j-1)_{k-1}P_k(i)p(i).
 \end{equation}
 \end{enumerate}
 \end{theorem}

 \begin{remark}
 \label{rem for P_k when finite S: delta>0 and Eq^[k]}
 \begin{enumerate}[itemsep=.5ex, wide, labelwidth=!, labelindent=0pt, label=\rm(\alph*), ref=\textcolor{black}{\alph*}]
 \item
 \label{rem for P_k when finite S: delta>0 and Eq^[k](a)}
 In \eqref{eq Rodiquez APP} when $k>M$, we have
 $\E\left[P_k^2(X)\right]=0$, since the polynomial $q^{[k]}$ vanishes identically on
 $S$. Thus, in the sequel, we study the polynomials $P_k$
 only when $k\le M$.
 \item
 \label{rem for P_k when finite S: delta>0 and Eq^[k](b)}
 Provided that $\E|X|^{2k}<\infty$ and $k\le M$, the
 quantities $1-j\delta$, $j=0,1,\ldots,2k-2$, are strictly
 positive. If $\delta\le0$, this is obvious. If $\delta>0$
 and $S$ is infinite, this follows from
 \Cref{lem moments and delta}; when $S$ is finite, it follows from
 \Cref{rem finite S: delta>0}. Thus, the quantity
 $\Pi_\delta^{[k]}(k-1)$
 is strictly positive.
 Also, since the polynomial $q^{[k]}$ is non-negative on $S$
 and $\Pr\left[q^{[k]}(X)>0\right]>0$, it follows that
 $0<\E\left[q^{[k]}(X)\right]<\infty$.
 \end{enumerate}
 \end{remark}

 For a non-negative integer $n$ such that $n\le M$ and
 $\E|X|^{2n}<\infty$,
 \Cref{rem for P_k when finite S: delta>0 and Eq^[k]}\eqref{rem for P_k when finite S: delta>0 and Eq^[k](b)}
 shows that we can define the standardized Rodrigues polynomials,
 \begin{equation}
 \label{eq phi_k}
 \phi_k(j)={[k!c_k(\delta)A_k]^{-1/2}}{P_k(j)},\quad k=0,1,\ldots,n.
 \end{equation}
 The set $\{\phi_k\}_{k=0}^{n}\subset L^2(\RR,X)$ is an
 orthonormal basis for all polynomials with degree at most
 $n$. Moreover, \eqref{eq lead P_k} shows that the leading
 coefficient is given by
 \begin{equation}
 \label{eq lead phi_k}
 \lead(\phi_k)\doteq d_k(\mu;q)=[{c_k(\delta)}/{(k!A_k)}]^{1/2}>0,\quad
 k=0,1,\ldots,n.
 \end{equation}

 Let $X$ be any rv of the CO family
 with $\E|X|^{2n}<\infty$,
 where $n$ is less than the cardinality of the support of
 $X$. It is well-known that we can always construct an
 orthonormal set of real polynomials up to order $n$. This
 construction is based on the first $2n$ moments of $X$ and
 is a by-product of the Gram-Schmidt orthonormalization
 process, applied to the linearly independent system
 $\{1,x,x^2,\ldots,x^n\}\subset L^2(\RR,X)$. The orthonormal
 polynomials are then uniquely defined, apart from the fact
 that we can multiply each polynomial by $\pm1$. It follows
 that the standardized Rodrigues polynomials $\phi_k$ of
 \eqref{eq phi_k} are the unique orthonormal polynomials
 that can be defined for a pmf $p\sim\CO$, provided that
 $\lead(\phi_k)>0$. Therefore, it is useful to express the
 $L^2$-norm of each $P_k$ in terms of the parameters
 $\delta$, $\beta$, $\gamma$ and $\mu$. This result is
 given by \eqref{eq orthogonality condition} and
 \Cref{prop q^[i](j)p(j)}\eqref{prop q^[i](j)p(j)(d)}.

 Consider the rvs $X_i$ with pmfs $p_i$ as defined in
 \Cref{prop q^[i](j)p(j)}.
 From \eqref{eq Rodiquez APP},
 the corresponding Rodrigues polynomials
 are given  by
 \begin{equation}
 \label{eq Rodiquez APP order i}
 P_{k,i}(j)=\frac{(-1)^k}{p_i(j)}\Delta^k\left[q_i^{[k]}(j-k)p_i(j-k)\right].
 \end{equation}
 Thus, the standardized Rodrigues polynomials,
 orthonormal with respect to the pmf of $X_i$,
 are given by
 \begin{equation}
 \label{eq phi_k,i}
 \phi_{k,i}(j)={\left[k!c_k(\delta_i)A_k(\mu_i;q_i)\right]^{-1/2}}{P_{k,i}(j)}.
 \end{equation}
 Note that for $i=1$, the rv $X_1$ is denoted by $X^*$
 ($p_1\equiv p^*$ etc.). Therefore, we may denote
 the polynomial $P_{k,1}$ by $P_k^*$ and the
 standardized polynomial $\phi_{k,1}$ by $\phi_k^*$.
 An important observation is that the
 forward difference of $\phi_k$ is
 scalar multiple of $\phi_{k-1}^*$. Specifically, we
 have the following lemma.

 \begin{lemma}
 \label{lem Delta phi_k}
 If $X\sim\COq=\CO$ and $\E|X|^{2n}<\infty$ for some $1\le n
 \le M$, then the polynomials $\phi_k$ of \eqref{eq
 phi_k} and $\phi_{k,1}\equiv\phi_k^*$ of \eqref{eq
 phi_k,i} are related through
 \begin{equation}
 \label{eq Delta phi_k}
 \begin{split}
 \Delta\phi_k(j)=v_{k-1}&\phi_{k-1}^*(j), \quad k=1,2,\ldots,n, \quad \textrm{where} \\
 & v_{k-1}=v_{k-1}(\mu;q)\doteq\left\{{k[1-(k-1)\delta]}/{A_1}\right\}^{1/2}.
 \end{split}
 \end{equation}
 \end{lemma}
 \begin{proof}
 First, we show that for $1\le m<k\le n$,
 $\E\left[\Delta\phi_k(X^*)\Delta\phi_m(X^*)\right]=0$.
 We have
 \begin{equation}
 \label{eq Delta phi_k(j) Delta phi_m(j) p(j)}
\!\! \scalebox{.95}{$[\Delta\phi_k(j)\Delta\phi_m(j)]q(j)p(j)
 =\Delta\{\phi_k(j)[\Delta\phi_m(j-1)]q(j-1)p(j-1)\}
 -\phi_k(j)\mathrm{pol}_m(j)p(j),$}
 \end{equation}
 where
 $\mathrm{pol}_m(j)\doteq\left[\Delta^2\phi_m(j-1)\right]q(j)+[\Delta\phi_m(j-1)](\mu-j)$
 is a polynomial with $\deg(\mathrm{pol}_m)\le{m}$.
 Summing \eqref{eq Delta
 phi_k(j) Delta phi_m(j) p(j)}
 for all $j\in\{\alpha,\ldots,\omega\}$, we observe the following:
 The lhs of the sum is
 $\E[\Delta\phi_k(X^*)\Delta\phi_m(X^*)]\E[q(X)]$. The
 first part of the rhs of the sum is
 $\phi_k(j)[\Delta\phi_m(j-1)]q(j-1)p(j-1)|_{\alpha}^{\omega+1}=0$
 (for finite $\alpha$ and $\omega$, this follows from
 $p(\alpha-1)=q(\omega)=0$; for infinite $\alpha$ and $\omega$,
 it follows from \Cref{prop lim_j->oo j^2ip(j)}). The
 second part of the rhs of the sum is $\E[\phi_k(X)\mathrm{pol}_m(X)]=0$,
 because
 $\phi_k$ is orthogonal to
 any polynomial of degree less than $k$.
 From the moment conditions, it is obvious that
 $\E[\Delta\phi_k(X^*)]^2<\infty$. Thus, it suffices
 to show that $\E[\Delta\phi_k(X^*)]^2>0$.
 The polynomial
 $\Delta\phi_k(x)$, $x\in\RR$, is not identically zero,
 since
 $\lead(\Delta\phi_k)=k\lead(\phi_k)>0$, and can not vanish
 identically on the support of
 $X^*$, since $\deg(\Delta\phi_k)=k-1$ is less than the
 cardinality of the support of $X^*$.
 Finally, since $\deg(\Delta\phi_k)=\deg\left(\phi^*_{k-1}\right)=k-1$,
 $k=1,\ldots,n$, the uniqueness of the orthogonal
 polynomial system implies that there exist constants
 $v_k\ne0$ such that $\Delta\phi_k=v_{k-1}\phi^*_{k-1}$.
 Equating the leading coefficients, we obtain
 $\lead(\Delta\phi_k)=v_{k-1}\lead(\phi^*_{k-1})$, that is,
  $v_{k-1}={\lead(\Delta\phi_k)}\left/{\lead\left(\phi^*_{k-1}\right)}\right.
        ={k\lead(\phi_k)}\left/{\lead\left(\phi^*_{k-1}\right)}\right.
        =k\{[(k-1)!c_k(\delta)A_{k-1}(\mu^*;q^*)]/[k!c_{k-1}(\delta^*)A_k]\}^{1/2}$;
 see \eqref{eq lead phi_k}. Moreover, one can easily see
 that
 $c_k(\delta)=[1-(k-1)\delta](1-2\delta)^{k-1}c_{k-1}(\delta^*)$
 and
 $A_k=(1-2\delta)^{k-1}A_1A_{k-1}(\mu^*;q^*)$.
 Thus,
 $v_{k-1}=\{k[1-(k-1)\delta]/A_1\}^{1/2}$.
 \end{proof}

 Applying now \Cref{lem Delta phi_k}, inductively it is easy
 to verify the following result.

 \begin{theorem}
 \label{theo Delta^m phi_k}
 Let $X\sim\COq=\CO$ and assume that  $\E|X|^{2n}<\infty$
 for some integer $n$ with $1\le n\le M$. Then,
 \begin{equation}
 \label{eq Delta^m phi_k}
 \begin{split}
 \Delta^m\phi_k(j&)=v^{(m)}_{k-m}\phi_{k-m,m}(j), \quad m=0,1,\ldots,n, \quad k=m,m+1,\ldots,n, \\
 &\textrm{with} \quad
 v^{(m)}_{k-m}=v^{(m)}_{k-m}(\mu;q)\doteq\left\{\left.{k!\Pi_\delta^{[m]}(k-1)}\right/[(k-m)!A_m]\right\}^{1/2},
 \end{split}
 \end{equation}
 where the polynomials $\phi_k$, $\phi_{k-m,m}$ are as given in
 \eqref{eq phi_k} and \eqref{eq phi_k,i}, respectively.
 \end{theorem}
 \begin{proof}
 The proof follows by induction on $m$. For $m=0$, the result
 is obvious, noting that $\phi_{k,0}=\phi_k$
 and $\nu_{k,0}^{(0)}=1$.
 For $m=1$, the result follows by \Cref{lem
 Delta phi_k}, since $\phi_{k,1}=\phi_k^*$
 and $\nu_{k,1}^{(1)}=\nu_k$.
 Assuming that it is true for
 $m-1\in\{0,1,\ldots,n-1\}$, we will show that it
 holds for $m$.
 By the assumption of induction,
 $\Delta^{m-1}\phi_k(j)=v^{(m-1)}_{k-m+1}\phi_{k-m+1,m-1}(j)$,
 and
 $v^{(m-1)}_{k-m+1}=\left\{\left.{k!\Pi_\delta^{[m-1]}(k-1)}\right/[(k-m+1)!A_{m-1}]\right\}^{1/2}$.
 Applying \Cref{lem Delta phi_k} for
 $X_{m-1}\sim\mathrm{CO}(\mu_{m-1};q_{m-1})$,
 $\Delta^{m}\phi_k(j)=\Delta\left[\Delta^{m-1}\phi_k(j)\right]
                     =v^{(m-1)}_{k-m+1}\Delta\left[\phi_{k-m+1,m-1}(j)\right]
                    =v^{(m-1)}_{k-m+1}v_{k-m}(\mu_{m-1};q_{m-1})\phi^*_{k-m,m-1}(j)
                     =v^{(m-1)}_{k-m+1}v_{k-m}(\mu_{m-1};q_{m-1})\phi_{k-m,m}(j)$,
 where
 $ v^{(m-1)}_{k-m+1}v_{k-m}(\mu_{m-1};q_{m-1})
 =\left(\left\{k!\left.\Pi_\delta^{[m-1]}(k-1)\right/[(k-m+1)!A_{m-1}]\right\}\right.
  \{(k-m+1)[1-(k-m)\delta_{m-1}]/[A_1(\mu_{m-1};q_{m-1})]\}\left.\vphantom{\Pi_\delta^{[m-1]}}\right)^{1/2}$;
 see \eqref{eq Delta phi_k}.
 Finally, it is easily shown that
 $A_1(\mu_{m-1};q_{m-1})={A_m}/\{[1-2(m-1)\delta]A_{m-1}\}$
 and $1-(k-m)\delta_{m-1}=[1-(k-2m-2)\delta]/[1-2(m-1)\delta]$.
 Thus,
 $v^{(m-1)}_{k-m+1}v_{k-m}(\mu_{m-1};q_{m-1})=v^{(m)}_{k-m}$,
 completing the proof.
 \end{proof}

 \section[$L^2$ completeness and expansions]{\bm{$L^2$} completeness and expansions}
 \label{section L^2 completeness and expansions}
 We now study the Fourier coefficients of a function
 regarding its expansion in the $L^2$ Hilbert space.
 First, we present
 the following basic result.
 \begin{theorem}[\textrm{\cite[Theorem 2.2]{APP2}}]
 \label{theo Fourier coefficients in CO}
 Suppose $X\sim\COq$ and that $\E|X|^{2k}<\infty$ for some $k\geq 1$.
 If $g$ is a function defined on $S$ with
 $\E\left[q^{[k]}(X)\left|\Delta^k g(X)\right|\right]<\infty$, then
 $\E|P_k(X)g(X)|<\infty$ and the following covariance
 identity holds:
 \begin{equation}
 \label{eq cov-ident generalize}
 \E[P_k(X)g(X)]=\E\left[q^{[k]}(X)\Delta^k g(X)\right].
 \end{equation}
 \end{theorem}
 Note that if the support $S$  has a finite upper endpoint,
 $\omega<\infty$, then $\Delta^k g(j)$, $j\in S$,
 may depend on some values $\{g(j), \ j\notin S\}$;
 however, only the values $\{j\colon j\in S, j\le\omega-k\}$
 are relevant to the rhs of the covariance identity \eqref{eq cov-ident generalize}.
 This is so because for $j>\omega-k$, the ascending power $q^{[k]}(j)$
 includes the factor $q(\omega)=0$.
 Thus, assuming any values for $g(j)$ when $j$ lies
 in the set $\{\omega+1,\omega+2,\ldots\}$,
 e.g., $g(j)=0$,
 $j=\omega+1,\omega+2,\ldots$, will not
 affect the covariance identity.
 For any
 function $g$ defined on $S$, the function
 $\Delta^k g$ has domain the set $S_k$;
 see \Cref{prop q^[i](j)p(j)}\eqref{prop q^[i](j)p(j)(a)}. Thus, the values $\Delta^k
 g(j)$, $j\in S\smallsetminus S_k$ (if exist),
 that appear in the formula, are immaterial.
 Note that if $S$ is finite and $k>M(X)$, then both polynomials
 $P_k$ and $q^{[k]}$ are identically zero
 on $S$, and the relation \eqref{eq cov-ident generalize} takes
 the trivial form $0=0$.

 It is important to note that the identity
 \eqref{eq cov-ident generalize}, combined
 with \eqref{eq phi_k},
 enables a convenient calculation
 of the Fourier coefficient $\alpha_k=\E[\phi_k(X)g(X)]$ of a
 function $g$. Specifically,
 \begin{equation}
 \label{eq E phi_k(X)g(X)}
 \alpha_k=\E[\phi_k(X)g(X)]={[k!c_k(\delta)A_k]^{-1/2}}{\E\left[q^{[k]}(X)\Delta^k g(X)\right]}.
 \end{equation}
 The rhs of \eqref{eq E phi_k(X)g(X)}
 shows that we do not need to know the polynomial $\phi_k$
 in order to calculate $\alpha_k$.

 We now shed
 some light on the interrelations
 between the spaces $L^2(\RR,X_i)$ and  $L^1(\RR,X_i)$.
 \begin{lemma}
 \label{lem L^2 and L^1}
 Let the rvs $X$ and $X^*$ be as in \Cref{lem q(j)p(j)}.
 Assume that the function $g$ is defined on the support of
 $X$.
 Then,
 \begin{enumerate}[itemsep=.5ex, wide, labelwidth=!, labelindent=0pt, label=\rm(\alph*), ref=\textcolor{black}{\alph*}]
 \item
 \label{lem L^2 and L^1(a)}
 $\Delta g \in L^2(\RR,X^*) \Rightarrow g \in L^2(\RR,X)$;
 \item
 \label{lem L^2 and L^1(b)}
 $\Delta g \in L^1(\RR,X^*) \Rightarrow g \in L^1(\RR,X)$.
 \end{enumerate}
 \end{lemma}
 \begin{proof}
 \eqref{lem L^2 and L^1(a)} For $|S|<\infty$, the result is obvious. Thus,
 assume that $|S|=\infty$ and consider a function $g$ such
 that  $\Delta g \in L^2(\RR,X^*)$.
 It suffices to show that for some $m\in\ZZ$,
 \[
 \sum_{\mathclap{j=m}}^{\infty}g^2(j)p(j)<\infty \ \ \textrm{when} \ \omega=\infty,
 \quad\textrm{and}\quad
 \sum_{\mathclap{j=-\infty}}^{m}g^2(j)p(j)<\infty \ \ \textrm{when} \
 \alpha=-\infty.
 \]
 For the first inequality, it suffices to show that
 $\Sigma_1(m)\doteq\sum_{j=m}^{\infty}[g(j)-g(m)]^2p(j)<\infty$.
 Let $m=\lfloor\mu\rfloor+1>\mu$. Then,
 $\Sigma_1(m)=\sum_{j=m}^{\infty}p(j)\left[\sum_{i=m}^{j-1}\Delta g(i)\right]^2
             \le\sum_{j=m}^{\infty}p(j)(j-m)\sum_{i=m}^{j-1}[\Delta g(i)]^2
            =\sum_{i=m}^{\infty}[\Delta g(i)]^2\sum_{j=i+1}^{\infty}(j-m)p(j)
             \le\sum_{i=m}^{\infty}[\Delta g(i)]^2\sum_{j=i+1}^{\infty}(j-\mu)p(j)$.
 Since $\sum_{j=i+1}^{\infty}(j-\mu)p(j)=q(i)p(i)$, we get
 $\Sigma_1(m)\le\sum_{i=m}^{\infty}[\Delta
 g(i)]^2q(i)p(i)\le\sum_{i\in\ZZ}[\Delta g(i)]^2q(i)p(i)=\E[q(X)]\E[\Delta
 g(X^*)]^2<\infty$.
 For the second inequality, we use the same arguments with
 $m=\lfloor\mu\rfloor\le\mu$.

 \eqref{lem L^2 and L^1(b)} Let $\Delta g \in L^1(\RR,X^*)$.
 Then, $\E[q(X)|\Delta g(X)|]=\E[q(X)]\E|\Delta g(X^*)|<\infty$.
 Applying \Cref{theo Fourier coefficients in CO} for $k=1$, and since
 $P_1(j)=j-\mu$, it follows that
 $\E|P_1(X)g(X)|=\sum_{j\in\ZZ}|(j-\mu)g(j)|p(j)$ is finite.
 Thus,
 $\sum_{j>\lfloor\mu\rfloor+1}|g(j)|p(j)\le\sum_{j>\lfloor\mu\rfloor+1}|(j-\mu)g(j)|p(j)<\infty$
 and
 $\sum_{j\le\lfloor\mu\rfloor-1}|g(j)|p(j)\le\sum_{j\le\lfloor\mu\rfloor-1}|(j-\mu)g(j)|p(j)<\infty$,
 completing the proof.
 \end{proof}

 \begin{corollary}
 \label{cor L^2 and L^1}
 Let the rvs $X$ and $X_i$, $i=0,1,\ldots,n$ be as in
 \Cref{prop q^[i](j)p(j)} and consider a function $g$
 defined on the support of $X$. Then:
 \begin{enumerate}[itemsep=.5ex, wide, labelwidth=!, labelindent=0pt, label=\rm(\alph*), ref=\textcolor{black}{\alph*}]
 \item
 \label{cor L^2 and L^1(a)}
 $\Delta^n g \in L^2(\RR,X_n) \Rightarrow \Delta^i
 g \in L^2(\RR,X_i)$ for every $i=0,1,\ldots,n$;
 \item
 \label{cor L^2 and L^1(b)}
 $\Delta^n g \in L^1(\RR,X_n) \Rightarrow \Delta^i g \in
 L^1(\RR,X_i)$ for every $i=0,1,\ldots,n$.
 \end{enumerate}
 \end{corollary}
 \begin{proof}
 Follows immediately by an application of \Cref{lem L^2 and L^1}.
 \end{proof}

 It is known (due to M. Riesz) that the real
 polynomials are dense in $L^2(\RR,X)$ whenever the
 probability measure of $X$ is determined by its moments;
 see \cite{Riesz,Akh}. An even simpler sufficient condition
 is when $X$ has a finite moment generating function at a
 neighborhood of zero, that is, when there exists $t_0>0$
 such that
 \begin{equation}
 \label{eq M_X(t)<infty}
 M_X(t)=\E\left(e^{tX}\right)<\infty,\quad t\in(-t_0,t_0);
 \end{equation}
 see \cite{APP2}, cf.\ \cite{BC}.

 Consider a rv $X$ in the CO family. If the support of $X$ is
 finite, then \eqref{eq M_X(t)<infty} holds, and obviously, the real
 polynomials are dense in the finite-dimensional space $L^2(\RR,X)$;
 in this case,
 $L^2(\RR,X)=\Span\left\{1,x,x^2,\ldots,x^M\right\}$, and the system of
 polynomials $\{\phi_k\}_{k=0}^{M}$ is an orthonormal basis
 of $L^2(\RR,X)$. When $X$ has infinite support, then there
 are two possibilities: If $\delta>0$, then $X$ does not have
 finite moments of any order,
 see \Cref{lem moments and delta}, and any real
 polynomial of $L^2(\RR,X)$ is of bounded degree; thus,
 only a finite number of orthonormal polynomials exist,
 and these polynomials cannot
 be dense in the infinite-dimensional space $L^2(\RR,X)$.
 If $\delta\le0$,
 then \eqref{eq M_X(t)<infty} holds, see \Cref{section classification}
 or \Cref{table probabilities}, so the real polynomials are dense in
 $L^2(\RR,X)$ and the system of polynomials
 $\{\phi_k\}_{k=0}^{\infty}$ is an orthonormal basis of this
 space. From the above observations, it is natural to
 define the following subclass of
 rvs of the CO system:
 \[
 \CX\doteq\{X\colon X\sim\CO \textrm{ for some } (\mu;\delta,\beta,\gamma),
 \textrm{ and }  \delta\le0 \textrm{ or } |S(X)|<\infty\}.
 \]
 \begin{remark}
 \label{rem X_i in C}
 Let $X\in\CX$. Then:
 \begin{enumerate}[itemsep=.5ex, wide, labelwidth=!, labelindent=0pt, label=\rm(\alph*), ref=\textcolor{black}{\alph*}]
 \item
 \label{rem X_i in C(a)}
 The set of polynomials $\{\phi_k\}_{k=0}^{M}$ ($M$ is
 finite or infinite) is an orthonormal basis of
 $L^2(\RR,X)$. Thus, any function $g\in L^2(\RR,X)$ can be
 expanded as
 \begin{equation}
 \label{eq g expand}
 g(j)\sim\sum_{\mathclap{k=0}}^{M}\alpha_k\phi_k(j),
 \end{equation}
 where $\alpha_k=\E[\phi_k(X)g(X)]$ are the Fourier coefficients of
 $g$. The series converges in the norm of $L^2(\RR,X)$; that
 is, $\E\left[g(X)-\sum_{k=0}^{M}\alpha_k\phi_k(X)\right]^2=0$ (when
 $M<\infty$) or
 $\E\left[g(X)-\sum_{k=0}^{N}\alpha_k\phi_k(X)\right]^2\to0$ as
 $N\to\infty$ (when $M=\infty$). Parseval's identity shows
 that
 \begin{equation}
 \label{eq Var g(X)=sum_k=0^M a_k^2}
 \Var[g(X)]=\sum_{\mathclap{k=1}}^{M}\alpha_k^2,\quad g\in L^2(\RR,X);
 \end{equation}
 \item
 \label{rem X_i in C(b)}
 For every $i=0,1,\ldots,M$, $X_i\in\CX$  (see \Cref{prop q^[i](j)p(j)}),
 and the corresponding results of
 \eqref{rem X_i in C(a)} hold for each $X_i$.
 \end{enumerate}
 \end{remark}

 One can apply $i$ times the forward difference operator in
 the series \eqref{eq g expand} to get, in view of
 \Cref{theo Delta^m phi_k}, the formal expansion
 \begin{equation}
 \label{eq Delta^i g expand}
 \Delta^i
 g(j)\sim\sum_{\mathclap{k=i}}^{M}\alpha_k\Delta^i\phi_k(j)
 =\sum_{\mathclap{k=i}}^{M}v_{k-i}^{(i)}(\mu;q)\alpha_k\phi_{k-i,i}(j),
 \end{equation}
 where $v_{k-i}^{(i)}(\mu;q)$ and $\{\phi_{k-i,i}(j)\}_{k=i}^M$ are given by
 \eqref{eq Delta^m phi_k} and \eqref{eq phi_k,i}, respectively.
 Now, if the expansion \eqref{eq Delta^i g expand} was indeed correct in the
 $L^2(\RR,X_i)$-sense,
 then the completeness of the system
 $\{\phi_{k,i}\}_{k=0}^{M_i}$ in $L^2(\RR,X_i)$
 would lead to
 the corresponding Parseval identity,
 \begin{equation}
 \label{eq E q^[i](Delta^i g)^2/E q^[i]}
 \frac{\E\left\{q^{[i]}(X)\left[\Delta^ig(X)\right]^2\right\}}{\E\left[q^{[i]}(X)\right]}
 =\E\left[\Delta^ig(X_i)\right]^2
 =\sum_{\mathclap{k=i}}^{M}\left[v_{k-i}^{(i)}(\mu;q)\right]^2\alpha_k^2.
 \end{equation}
 Finally, from \eqref{eq Delta^m phi_k}, we have
 $\left[v_{k-i}^{(i)}(\mu;q)\right]^2=\left.{k!\Pi_\delta^{[i]}(k-1)}\right/\left\{(k-i)!\E\left[q^{[i]}(X)\right]\right\}$.
 A combination of the last equation with \eqref{eq E
 q^[i](Delta^i g)^2/E q^[i]} yields the important
 identity
 \begin{equation}
 \label{eq E q^[i](Delta^i g)^2}
 \E\left\{q^{[i]}(X)\left[\Delta^ig(X)\right]^2\right\}
 =\sum_{\mathclap{k=i}}^{M}\frac{k!\Pi_\delta^{[i]}(k-1)}{(k-i)!}\alpha_k^2.
 \end{equation}
 This should be correct for all $g$ such that $\Delta^i g\in
 L^2(\RR,X_i)$, provided that expansion \eqref{eq g expand}
 is valid. We shall show that this is indeed the case.
 The $L^2$
 convergence
 of $\sum_{k=0}^{N}\alpha_k\phi_k(X)$ to
 $g(X)$ implies that
 $g(X)=\sum_{k=0}^{M}\alpha_k\phi_k(X)$
 with probability 1, that is,
 $g(j)=\sum_{k=0}^{M}\alpha_k\phi_k(j)$ for all $j\in S(X)$.
 Therefore, $\Delta^i
 g(j)=\sum_{k=0}^{M}\alpha_k\Delta^i\phi_k(j)
 =\sum_{k=i}^{M}v_{k-i}^{(i)}(\mu;q)\alpha_k\phi_{k-i,i}(j)$
 for all $j\in S(X_i)$.

 However, the same result can be derived by an alternative
 technique, similar to the one given in \cite{AP1}. In fact,
 we shall show more, namely, that an initial segment of the
 Fourier coefficients for the $i$th difference of $g$,
 suggested by \eqref{eq Delta^i g expand}, can be derived
 for any $X\sim\CO$ having a sufficient number of moments.
 This result holds even if $\delta>0$ and $|S|=\infty$. We
 present this technique
 since \Cref{lem extended CO equation} and \Cref{theo Fourier alternative}
 may be of interest on their own right.

 \begin{lemma}
 \label{lem discrete drunken}
 Consider a non-negative sequence $\{a_i\}_{i\in\ZZ}$ and assume
 that there is a positive integer $n$ such that $\sum_{i\in\ZZ}|i|^n a_i$ is
 finite. For each $k\in\{0,1,\ldots,n\}$, we define
 the sequence $\{b_{j;k}\}_{j\in\ZZ}$ by the relation
 $b_{j;k}\doteq\sum_{i\ge{j}}[j-i]_k a_i$.
 Then:
 \begin{enumerate}[itemsep=.5ex, wide, labelwidth=!, labelindent=0pt, label=\rm(\alph*), ref=\textcolor{black}{\alph*}]
 \item
 \label{lem discrete drunken(a)}
 For every $k\in\{1,2,\ldots,n\}$,
 $\Delta b_{j;k}=k b_{j+1;k-1}$, where the forward
 difference is taken with respect to the index $j$;
 \item
 \label{lem discrete drunken(b)}
 $\Delta^r b_{j;n}=(n)_r b_{j+r;n-r}$ for each
 $r\in\{1,2,\ldots,n\}$. In particular, for
 $r=n$,
 \[
 \Delta^n b_{j;n}=n! b_{j+n;0}=n!\sum_{\mathclap{i\ge{j+n}}}a_i.
 \]
 \end{enumerate}
 \end{lemma}
 \begin{proof}
 \eqref{lem discrete drunken(a)}
 $\Delta b_{j;k}=\sum_{i\ge{j+1}}[j+1-i]_k
 a_i-\sum_{i\ge{j}}[j-i]_k a_i=\sum_{i\ge{j+1}}\Delta[j-i]_k
 a_i-[0]_k$. Since $[0]_k=0$ ($k>0$) and
 $\Delta[j-i]_k=k[j+1-i]_{k-1}$, the desired result follows.

 \eqref{lem discrete drunken(b)}
 It follows easily by applying \eqref{lem discrete drunken(a)} $r$ times inductively.
 \end{proof}

 \begin{lemma}
 \label{lem extended CO equation}
 Let $X\sim\COq=\CO$ and consider a positive integer $k\le
 M$.
 Then, provided that $\E|X|^{2k-1}$ is finite,
 \[
 q(j)p(j)\Delta
 P_k(j)=-\lambda_k(\delta)\sum_{i\le{j}}P_k(i)p(i)=\lambda_k(\delta)\sum_{i>j}P_k(i)p(i),
 \]
 where $\lambda_k(\delta)\doteq k[1-(k-1)\delta]$ and $P_k$
 is the
 orthogonal polynomial given by
 \eqref{eq Rodiquez APP}. If, in addition, $\E|X|^{2k}$ is
 finite, then for the standardized
 polynomial $\phi_k=\left\{\E\left[P_k^2(X)\right]\right\}^{-1/2}P_k$,
 we have
 \begin{equation}
 \label{eq extended CO standardized}
 q(j)p(j)\Delta\phi_k(j)
 =-\lambda_k(\delta)\sum_{i\le{j}}\phi_k(i)p(i)
 =\lambda_k(\delta)\sum_{i>j}\phi_k(i)p(i).
 \end{equation}
 \end{lemma}
 \begin{proof}
 Since $(x)_n=(-1)^n[-x]_n$, applying \eqref{eq Rodrigues inversion formula} (replacing $j$ by $j-(k-1)$),
 \begin{equation}
 \label{eq Rodrigues inversion formula2}
 q^{[k]}(j-(k-1))p(j-(k-1))
 =\frac{(-1)^{k-1}}{(k-1)!}\sum_{\mathclap{\quad i\ge{j-k+2}}}[j-i-k-2]_{k-1}P_k(i)p(i).
 \end{equation}
 The lhs of \eqref{eq Rodrigues inversion formula2}
 can be written as
 $
                              (1-2\delta)^{k-1}q_1^{[k-1]}(j-(k-1))p_1(j-(k-1))
                            \E[q(X)]$.
 Applying the operator $\Delta^{k-1}$
 and using \eqref{eq Rodiquez APP order i},
 we obtain
 $(-1)^{k-1}(1-2\delta)^{k-1}p_1(j)P_{k-1,1}(j)\E[q(X)]
 =(-1)^{k-1}(1-2\delta)^{k-1}q(j)p(j)P_{k-1,1}(j)$.
 As in \Cref{lem Delta phi_k},
 we find that
 $\Delta
 P_k(j)=B_{k-1}P_{k-1,1}(j)$, where
 $B_{k-1}
 ={\lead(\Delta P_k)}/{\lead(P_{k-1,1})}
 ={k\lead(P_k)}/{\lead(P_{k-1,1})}
 ={kc_k(\delta)}/{c_{k-1}(\delta_1)}
 =k[1-(k-1)\delta](1-2\delta)^{k-1}$.
 Therefore, an application of
 the operator $\Delta^{k-1}$
 to the lhs of \eqref{eq Rodrigues inversion formula2}
 produces the quantity
 ${(-1)^{k-1}}{\lambda_k^{-1}(\delta)}q(j)p(j)\Delta P_{k}(j)$.
 Applying the operator $\Delta^{k-1}$
 to the rhs of \eqref{eq Rodrigues inversion formula2}
 and using \Cref{lem discrete drunken}, we
 arrive at the quantity
 $(-1)^{k-1}\sum_{i>j}P_k(i)p(i)$,
 and the result follows from the fact that
 the last two quantities must be equal to each other.
 Finally, since $\E[P_k(X)]=0$  (because $k\ge1$), we conclude
 that
 $(-1)^{k-1}\sum_{i>j}P_k(i)p(i)=(-1)^{k}\sum_{i\le{j}}P_k(i)p(i)$.
 \end{proof}

 \begin{lemma}
 \label{lem Fourier alternative}
 Let the rvs $X$ and $X^*$ be as in \Cref{lem q(j)p(j)}, and assume that for some integer $k$ with $1\le k \le M$,
 $\E|X|^{{\max\{2k,3\}}}<\infty$. Then, for any function $g$ with
 $\Delta g\in L^2(\RR,X^*)$, we have the identity
 \begin{equation}
 \label{eq Fourier alternative 1}
 \E\left[\phi_{k-1}^*(X^*)\Delta g(X^*)\right]=v_{k-1}\E[\phi_{k}(X)g(X)],
 \end{equation}
 where $\phi_k$, $\phi_{k,1}\equiv\phi_k^*$ and
 $v_{k-1}=v_{k-1}(\mu;q)$ are as given in \eqref{eq
 phi_k}, \eqref{eq phi_k,i} and \eqref{eq Delta
 phi_k}, respectively.
 \end{lemma}
 \begin{proof}
 By an application of Cauchy-Schwarz inequality, we get
 $\E^2|\phi_{k-1}^*(X^*)\Delta g(X^*)|
 \le\E\left[\phi_{k-1}^*(X^*)\right]^2\E[\Delta g(X^*)]^2
 =\E[\Delta g(X^*)]^2<\infty$.
 From \Cref{cor L^2 and L^1}, it follows
 that $g\in L^2(\RR,X)$,
 and similarly, $\E|\phi_{k}(X)g(X)|<\infty$.
 Since $\E[\phi_k(X)]=0$, $\phi_k$ must change its sign
 in the support
 of $X$. Thus, $\phi_k$ has real roots, say $\rho_1<\cdots<\rho_m$,
 that lie in the interval $[\alpha,\omega]$. Fix now an
 integer $\rho\in\{[\rho_1],\ldots,[\rho_m]\}\subset S$.
 Then,
 $\E[q(X)] \E[\phi_{k-1}^*(X^*)\Delta g(X^*)]
    =\sum_{j=\alpha}^{\omega-1}\Delta g(j)q(j)p(j)\phi_{k-1}^*(j)
    =v_{k-1}^{-1}\sum_{j=\alpha}^{\omega-1}\Delta g(j)q(j)p(j)\Delta\phi_{k}(j)
    =-\lambda_k(\delta)v_{k-1}^{-1}\sum_{j=\alpha}^{\rho-1}\Delta g(j)\sum_{i=\alpha}^{j}p(i)\phi_k(i)
     +\lambda_k(\delta)v_{k-1}^{-1}\sum_{j=\rho}^{\omega-1}\Delta g(j)\sum_{i=j+1}^{\omega}p(i)\phi_k(i)$.
 Observing that $\lambda_k(\delta)v_{k-1}^{-1}=v_{k-1}\E[q(X)]$,
 the preceding equation can be rewritten as
 \[
 \E[\phi_{k-1}^*(X^*)\Delta g(X^*)]=v_{k-1}(\varSigma_2-\varSigma_1),
 \quad\textrm{where}
 \]
 \begin{equation}
 \label{eq Sigma1: Sigma2}
 \varSigma_1\doteq\sum_{\mathclap{j=\alpha}}^{\mathclap{\rho-1}}\Delta g(j)\sum_{\mathclap{i=\alpha}}^{j}p(i)\phi_k(i),
 \quad
 \varSigma_2\doteq\sum_{\mathclap{j=\rho}}^{\mathclap{\omega-1}}\Delta g(j)\sum_{\mathclap{i=j+1}}^{\omega}p(i)\phi_k(i).
 \end{equation}
 Now, we wish to change the order of summation to both sums
 $\varSigma_1$ and $\varSigma_2$. To this end, for
 $\varSigma_2$, it suffices to show that
 \begin{equation}
 \label{eq Sigma^*2}
 \varSigma^*_2\doteq\sum_{\mathclap{j=\rho}}^{\mathclap{\omega-1}}|\Delta
 g(j)|\sum_{\mathclap{i=j+1}}^{\omega}p(i)|\phi_k(i)|<\infty.
 \end{equation}
 Similarly, for $\varSigma_1$, it suffices to show that
 $\varSigma^*_1\doteq\sum_{j=\alpha}^{\rho-1}|\Delta
 g(j)|\sum_{i=\alpha}^{j}p(i)|\phi_k(i)|<\infty$. Note that,
 obviously, if $\alpha>-\infty$, then $\varSigma^*_1<\infty$
 and if $\omega<\infty$, then $\varSigma^*_2<\infty$. We now
 proceed to verify \eqref{eq Sigma^*2} when $\omega=\infty$.
 Write $\varSigma^*_2=\varSigma^*_{21}+\varSigma^*_{22}$,
 where
 $\varSigma^*_{21}\doteq\sum_{j=\rho}^{[\rho_m]}|\Delta g(j)|\sum_{i=j+1}^{\infty}p(i)|\phi_k(i)|$,
 and
 $\varSigma^*_{22}\doteq\sum_{j=[\rho_m]+1}^{\infty}|\Delta g(j)|\sum_{i=j+1}^{\infty}p(i)|\phi_k(i)|$.
 Since $\E|X|^{k}<\infty$,
 $\sum_{i=j+1}^{\infty}p(i)|\phi_k(i)|<\infty$ for each
 $j=\rho,\ldots,[\rho_m]$ and thus,
 $\varSigma^*_{21}<\infty$, being a finite sum of finite terms.
 On the other hand, since the polynomial $\phi_k$ does not
 change its sign in the set $\{[\rho_m]+1,[\rho_m]+2,\ldots\}$,
 we can define the constant
 $c\doteq\sign \phi_k(j) \in \{-1,1\}, \quad j\in\{[\rho_m]+1,[\rho_m]+2,\ldots\}$.
 Then, $c\phi_k(j)=|\phi_k(j)|$ holds for all
 $j\in\{[\rho_m]+1,[\rho_m]+2,\ldots\}$ and from
 \eqref{eq  extended CO standardized}, we get
 $\varSigma^*_{22}=c\sum_{j=[\rho_m]+1}^{\infty}|\Delta g(j)|\sum_{i=j+1}^{\infty}p(i)\phi_k(i)
                  =c\lambda_k^{-1}(\delta)\sum_{j=[\rho_m]+1}^{\infty}|\Delta g(j)|q(j)p(j)\Delta\phi_k(j)
                  \le\lambda_k^{-1}(\delta)\sum_{j=[\rho_m]+1}^{\infty}|\Delta g(j)|q(j)p(j)|\Delta\phi_k(j)|
                  \le\lambda_k^{-1}(\delta)\sum_{j=\alpha}^{\infty}|\Delta g(j)|q(j)p(j)|\Delta\phi_k(j)|
                  =v_{k-1}\lambda_k^{-1}(\delta)\E[q(X)]\sum_{j=\alpha}^{\infty}\left|\Delta g(j)\phi_{k-1}^*(j)\right|p^*(j)
                  =v_{k-1}^{-1}\E\left|\phi_{k-1}^*(X^*)\Delta g(X^*)\right|<\infty$.
 Therefore, \eqref{eq Sigma^*2} follows for both cases
 ($\omega<\infty$ or $\omega=\infty$). If $\alpha=-\infty$,
 using similar arguments it can be shown that
 $\varSigma^*_1<\infty$.
 Thus, we can indeed interchange the order of summation to
 both sums $\varSigma_1$ and $\varSigma_2$ of \eqref{eq
 Sigma1: Sigma2}. It follows that
 $\varSigma_1=\sum_{i=\alpha}^{\rho-1}p(i)\phi_k(i)\sum_{j=i}^{\rho-1}\Delta g(j)
            =g(\rho)\sum_{i=\alpha}^{\rho-1}p(i)\phi_k(i)-\sum_{i=\alpha}^{\rho-1}g(i)p(i)\phi_k(i)$ and
$\varSigma_2=\sum_{i=\rho+1}^{\omega}p(i)\phi_k(i)\sum_{j=\rho}^{i-1}\Delta g(j)
            =\sum_{i=\rho+1}^{\omega}g(i)p(i)\phi_k(i)-g(\rho)\sum_{i=\rho+1}^{\omega}p(i)\phi_k(i)
            =\sum_{i=\rho}^{\omega}g(i)p(i)\phi_k(i)-g(\rho)\sum_{i=\rho}^{\omega}p(i)\phi_k(i)$.
 Taking into account the fact that
 $\sum_{\alpha}^{\omega}p(i)\phi_k(i)=\E[\phi_k(X)]=0$, we get
$\varSigma_2-\varSigma_1=\sum_{\alpha}^{\omega}g(i)p(i)\phi_k(i)-g(\rho)\sum_{\alpha}^{\omega}p(i)\phi_k(i)=\E[\phi_k(X)g(X)]$,
 which completes the proof of the lemma.
 \end{proof}

 \begin{theorem}
 \label{theo Fourier alternative}
 Let $X\sim\COq=\CO$ and fix an integer $k$ with
 $1\le k \le M$. Assume
 that $\E|X|^{2k+1}<\infty$ and
 consider the rvs $X_i$, $i=0,1,\ldots,k$,
 as in \Cref{prop q^[i](j)p(j)}.
 Then:
 \begin{enumerate}[itemsep=.5ex, wide, labelwidth=!, labelindent=0pt, label=\rm(\alph*), ref=\textcolor{black}{\alph*}]
 \item
 \label{theo Fourier alternative(a)}
 The  Fourier coefficients satisfy
 the relation
 \begin{equation}
 \label{eq Fourier alternative i}
 \E\left[\phi_{k-i,i}(X_i)\Delta^ig(X_i)\right]=v_{k-i}^{(i)}\E[\phi_{k}(X)g(X)], \quad i=0,1,\ldots,k,
 \end{equation}
 where $\phi_k$, $\phi_{k,i}$ and
 $v_{k-i}^{(i)}=v_{k-i}^{(i)}(\mu;q)$ are as given in
 \eqref{eq phi_k}, \eqref{eq phi_k,i} and
 \eqref{eq Delta^m phi_k}, respectively;
 \item
 \label{theo Fourier alternative(b)}
 If, in addition, $X\in\CX$ and $\Delta^n g\in L^2(\RR,X_n)$
 for some fixed integer $n$ with $1\le n\le M$,
 then \eqref{eq E q^[i](Delta^i g)^2}
 holds for all $i=0,1,\ldots,n$.
 \end{enumerate}
 \end{theorem}
 \begin{proof}
 \eqref{theo Fourier alternative(a)} By \Cref{cor L^2 and L^1}, $\Delta^i g\in
 L^2(\RR,X_i)$ for all $i=0,1,\ldots,k$. For $i=0$,
 \eqref{eq Fourier
 alternative i} is obvious and for $i=1$, it
 follows from \Cref{lem Fourier alternative}. Assume that
 it is true for $i-1\in\{0,\ldots,k-1\}$, that is,
 $\E\left[\phi_{k-i+1,i-1}(X_{i-1})\Delta^{i-1}g(X_{i-1})\right]=v_{k-i+1}^{(i-1)}\E[\phi_{k}(X)g(X)]$.
 Observe that the assumptions of
 \Cref{lem Fourier alternative} are satisfied
 for the rv $X_{i-1}$, the integer $k-i+1$
 and the function $\Delta^{i-1} g$.
 Using \eqref{eq Fourier alternative 1},
 $\E\left[\phi_{k-i,i}(X_{i})\Delta^{i}g(X_{i})\right]=\E\left[\phi_{k-i,i}(X_{i})\Delta\left(\Delta^{i-1}g(X_{i})\right)\right]
                                         =v_{k-i}(\mu_{i-1};q_{i-1})\E\left[\phi_{k-i+1,i-1}(X_{i-1})\Delta^{i-1}g(X_{i-1})\right]$.
 Thus, we get
 $\E\left[\phi_{k-i,i}(X_{i})\Delta^{i}g(X_{i})\right]=v_{k-i}(\mu_{i-1};q_{i-1})v_{k-i+1}^{(i-1)}\E[\phi_{k}(X)g(X)]$.
 Finally,
 $v_{k-i}(\mu_{i-1};q_{i-1})=\{(k-i+1)[1-(k-i)\delta_{i-1}]/{A_1(\mu_{i-1};q_{i-1})}\}^{1/2}$,
 where
 $A_1(\mu_{i-1};q_{i-1})={A_i}/\{[1-2(i-1)\delta]A_{i-1}\}$
 and
 $1-(k-i)\delta_{i-1}=[1-(k+i-2)\delta]/[1-2(i-1)\delta]$.
 Hence,
 $v_{k-i}(\mu_{i-1};q_{i-1})=\{(k-i+1)[1-(k+i-2)\delta]A_{i-1}/A_{i}\}^{1/2}$
 and a straightforward calculation gives
 $v_{k-i}(\mu_{i-1};q_{i-1})v_{k-i+1}^{(i-1)}=v_{k-i}^{(i)}$.

 \eqref{theo Fourier alternative(b)}
 Since $X\in\CX$, we have that $X_i\in\CX$ and the set
 of polynomials $\{\phi_{k,i}\}_{k=0}^{M_i}$ (where
 $M_i=M(X_i)=M-i$) is an orthonormal basis of
 $L^2(\RR,X_i)$; see \Cref{rem X_i in C}\eqref{rem X_i in C(b)}. Moreover,
 $\Delta^i g\in L^2(\RR,X_i)$. Thus, by Parseval's identity, it follows that
 $\E\left[\Delta^ig(X_i)\right]^2=\sum_{k=0}^{M_i}\alpha_{k,i}^2
 =\sum_{k=i}^{M}\alpha_{k-i,i}^2$,
 where $\alpha_{k,i}\doteq\E\left[\phi_{k,i}(X_i)\Delta^ig(X_i)\right]$ (with
 $\alpha_{k,0}=\alpha_{k}$) is the Fourier coefficient of $\Delta^ig$
 with respect to $\phi_{k,i}$. Using \eqref{eq Fourier
 alternative i},
 $\alpha_{k-i,i}^2=\E^2\left[\phi_{k-i,i}(X_i)\Delta^ig(X_i)\right]
 =\left[v_{k-i}^{(i)}\right]^2\E^2[\phi_{k}(X)g(X)]
 =\left[v_{k-i}^{(i)}\right]^2\alpha_k^2$,
 which verifies \eqref{eq E q^[i](Delta^i g)^2/E q^[i]} and the
 proof is complete.
 \end{proof}

 \section{Applications to variance bounds}
 \label{section applications to variance bounds}
 We now use the results of \Cref{section L^2 completeness and expansions} to present a wide class of
 variance bounds for a function $g$ of a rv $X$
 in the CO family.

 Let $X$ be any rv in the CO family and consider
 two non-negative integers $m,n\le M$ such that
 $\E|X|^{2\ell}<\infty$, where
 $\ell=\max\{m,n\}$. We denote by $\CH^{m,n}(X)$ the
 class of functions
 $g\colon S\to\RR$ ($S=S(X)$ is the support of $X$)
 satisfying the restrictions
 \[
 \E\left\{q^{[n]}(X)\left[\Delta^ng(X)\right]^2\right\}<\infty
 \quad\textrm{and}\quad
 \E\left[q^{[m]}(X)|\Delta^m g(X)|\right]<\infty.
 \]
 From \Cref{cor L^2 and L^1} and
 the fact that
 $\E^2\left[q^{[i]}(X)|\Delta^ig(X)|\right]\le\E\left\{q^{[i]}(X)\left[\Delta^ig(X)\right]^2\right\}\times\E\left[q^{[i]}(X)\right]$
 for all $i=0,1,\ldots,n$,
 we conclude the following:
 \[
 \mbox{
 If $m\le n$ and if $\E
 \left\{q^{[n]}(X)\left[\Delta^ng(X)\right]^2\right\}<\infty$, then
 $\E\left[q^{[m]}(X)|\Delta^mg(X)|\right]<\infty$.
 }
 \]
 Note that \Cref{cor L^2 and L^1}
 requires $\E|X|^{2\ell+1}<\infty$,
 but this assumption is needed only for the existence
 of the pmf $p_\ell$; thus, for the validity
 of the above observation, it is
 sufficient that $\E|X|^{2\ell}<\infty$. It follows that
 $\CH^{0,n}=\CH^{1,n}=\cdots=\CH^{n,n}$ [of course,
 $\CH^{0,0}(X)=L^2(\RR,X)$].

 Furthermore, when $M=\infty$ and $X$ has finite moments of
 any order (that is, $\delta\le0$), we shall denote by
 $\CH^{\infty,n}(X)$ and
 $\CH^{\infty}(X)$
 the classes $\bigcap_{m=0}^{\infty}\CH^{m,n}(X)=\bigcap_{m=n+1}^{\infty}\CH^{m,n}(X)$
 and  $\bigcap_{n=0}^{\infty}\CH^{\infty,n}(X)$, respectively.
 That is,
 \[
 \begin{split}
 &\CH^{\infty,n}(X)=\left\{g\colon \E\left\{{q^{[n]}(X)}
 \left[\Delta^ng(X)\right]^2\right\}<\infty \textrm{ and }
 \E\left[{q^{[m]}(X)}|\Delta^m g(X)|\right]<\infty \ \forall m>n\right\},\\
 &\CH^{\infty}(X)=\left\{g\colon \E\left\{{q^{[n]}(X)}\left[\Delta^ng(X)\right]^2\right\}<\infty \ \forall n\in\N\right\}.
 \end{split}
 \]
 Note that, by definition,
 $\CH^{m,\infty}(X)\doteq\bigcap_{n=0}^\infty\CH^{m,n}
 \equiv\CH^{\infty}(X)$
 for arbitrary fixed $m$.

 From \Cref{cor L^2 and L^1},  we conclude that the
 (finite or infinite) sequence $\CH^{m,n}(X)$ is decreasing
 in both $m$ and $n$. In particular, if all moments of $X$
 exist, then
 \[
 \begin{array}{c@{\hspace{.3ex}}c@{\hspace{.3ex}}c@{\hspace{.3ex}}c@{\hspace{.3ex}}c@{\hspace{.3ex}}c@{\hspace{.3ex}}c@{\hspace{.3ex}}c@{\hspace{.3ex}}c@{\hspace{.3ex}}c}
  L^2(\RR,X)\equiv & \CH^{0,0}(X) &           &              &           &              &           &        &           &                 \\
 [-.5ex]
                  &  \invsubset  &           &              &           &              &           &        &           &                 \\
                  & \CH^{1,0}(X) & \supseteq & \CH^{1,1}(X) &           &              &           &        &           &                 \\
 [-.5ex]
                  &  \invsubset  &           & \invsubset   &           &              &           &        &           &                 \\
                  & \CH^{2,0}(X) & \supseteq & \CH^{2,1}(X) & \supseteq & \CH^{2,2}(X) &           &        &           &                 \\
[-.5ex]
                  &  \invsubset  &           & \invsubset   &           & \invsubset   &           &        &           &                 \\
[-1.5ex]
                  &  \vdots      &           & \vdots       &           & \vdots       &           &        &           &                 \\
[-.5ex]
                  &  \invsubset  &           & \invsubset   &           & \invsubset   &           &        &           &                 \\
                  & \CH^{M,0}(X) & \supseteq & \CH^{M,1}(X) & \supseteq & \CH^{M,2}(X) & \supseteq & \cdots & \supseteq & \CH^{M,M}(X).\\
\end{array}
\]

 \Cref{eq cov-ident generalize,eq E q^[i](Delta^i g)^2}
 are almost identical with those given in \cite[Eq.s (2.3) and (2.2)]{Afe},
 for the continuous case.
 Therefore, using similar arguments,
 the next theorem holds;
 cf. \cite[Theorem 2.1]{Afe}.

 \begin{theorem}
 \label{theo variance bounds}
 Let $X\in\CX$, and fix two non-negative integers $m,n$
 with  $1\leq m+n\le M$. Assume that the function
 $g\in\CH^{m,n}(X)$. Consider the quantity
 \begin{equation}
 \label{eq S_m:n}
 S_{m,n}(g)=\sum_{i=1}^{m}\kappa_i\E^2\left[q^{[i]}(X)\Delta^ig(X)\right]
 +\sum_{i=1}^{n}(-1)^{i-1}\nu_i\E\left\{{q^{[i]}(X)\left[\Delta^ig(X)\right]^2}\right\},
 \end{equation}
 where
 \[
 \kappa_i\doteq\frac{{m\choose{i}}\Pi_\delta^{[n]}(m+i)}{(m+n)_i\Pi_\delta^{[i]}(i-1)\Pi_\delta^{[n]}(m)\E\left[{q^{[i]}(X)}\right]}
 \quad \textrm{and} \quad \nu_i\doteq\frac{{n\choose{i}}}{(m+n)_i\Pi_\delta^{[i]}(m)}
 \]
 are strictly positive constants
 (depending only on $m,n$ and $X$),
 and an empty sum (when $m=0$ or $n=0$)
 should be treated as zero.
 Then, the following inequality holds:
 \[
 (-1)^n\{\Var[g(X)]-S_{m,n}(g)\}\ge0.
 \]
 Moreover, $S_{m,n}(g)$ becomes equal to $\Var[g(X)]$ if and
 only if $g$ is identically equal to a polynomial of degree at
 most $m+n$ on the support of $X$, that is, if and only if there
 exists a polynomial  $H_{m+n}$ of degree at most $m+n$ such that
 $\Pr[g(X)=H_{m+n}(X)]=1$.
 \end{theorem}
 \begin{proof}
 Let $\alpha_k=\E[\phi_k(X)g(X)]$ be the Fourier coefficients
 of $g$.
 From \eqref{eq E q^[i](Delta^i g)^2} and \eqref{eq cov-ident generalize},
 we get,  as in \cite{Afe}, that
 $(-1)^n\{\Var[g(X)]-S_{m,n}(g)\}=R_{m,n}(g)$, where
 \begin{equation}
 \label{eq R_m:n}
 R_{m,n}(g)=\sum_{\mathclap{k=m+n+1}}^{M}r_{k;m,n}(\delta)\alpha_k^2
 \doteq\sum_{\mathclap{k=m+n+1}}^{M}\frac{(k-m-1)_n\Pi_\delta^{[n]}(m+k)}
 {(m+n)_n\Pi_\delta^{[n]}(m)}\alpha_k^2.
 \end{equation}
 If $\delta\le0$, $\Pi_\delta^{[n]}(m+k)>0$ and $\Pi_\delta^{[n]}(m)>0$
 because $1-j\delta\ge1$ for all $j\in\NN$,
 while if $\delta>0$, the same follows by \Cref{rem finite S: delta>0}.
 Therefore, the residual $R_{m,n}(g)$ in \eqref{eq R_m:n} is
 non-negative, and it is equal to zero if and only if $\alpha_k=0$ for
 all $k>m+n$, i.e., if and only if the function $g\colon S(X)\to\RR$
 is a polynomial of degree at most $m+n$. Note that if
 $m+n=M$ (in the case where $M$ is finite), the sum in
 \eqref{eq R_m:n} is empty and it is treated as zero.
 \end{proof}

 \begin{example}
 \label{exm Poisson}
 Suppose $X\sim\mathrm{Poisson}(\lambda)$ and consider a function $g\colon \N\to\R$.
 \Cref{theo variance bounds} produces the inequality
 $(-1)^n\{\Var[g(X)]-S_{m,n}(g)\}\ge0$,
 where
 \[
 \begin{split}
 S_{m,n}(g)=\sum_{i=1}^m \frac{\lambda^i}{i!} \frac{{m\choose i}}{{m+n \choose
 i}} \E^2\left[\Delta^i g(X)\right] +
 \sum_{i=1}^n (-1)^{i-1}\frac{\lambda^i}{i!} \frac{{n\choose i}}{{m+n \choose
 i}} \E \left[\Delta^i g(X) \right]^2, \\
 n,m=0,1,\ldots, \quad n+m>0,
 \end{split}
 \]
 provided $\E[\Delta^n g(X)]^2<\infty$ and
 $\E \left|\Delta^m g(X) \right|<\infty$ (of course, if $m\leq n$,
 the second restriction is implied by the first one). The equality
 holds if and only if $g\colon \N\to\R$ is a polynomial of degree at most
 $n+m$. For $n=m=1$, we get
 \eqref{eq new bound for Poisson m=n=1}.
 \end{example}

 \begin{remark}
 \label{rem S_m:n with n fixed}
 \begin{enumerate}[itemsep=.5ex, wide, labelwidth=!, labelindent=0pt, label=\rm(\alph*), ref=\textcolor{black}{\alph*}]
 \item
 \label{rem S_m:n with n fixed(a)}
 For fixed $n$ and for any function
 $g\in\CH^{\tilde{m},n}(X)$, where $\tilde{m}$ can
 be finite or infinite, the variance bounds
 $\{S_{m,n}(g)\}_{m=0}^{\tilde{m}}$ are of the same
 kind, i.e., upper bounds when $n$ is odd and lower bounds when
 $n$ is even;
 \item
 \label{rem S_m:n with n fixed(b)}
 The bounds $\{S_{m,n}(g)\}_{m=0}^{n}$ require the same
 condition on $g$, i.e., $g\in\CH^{n,n}(X)$.
 \end{enumerate}
 \end{remark}

 \begin{remark}
 \label{rem showed and existings}
 \begin{enumerate}[itemsep=.5ex, wide, labelwidth=!, labelindent=0pt, label=\rm(\alph*), ref=\textcolor{black}{\alph*}]
 \item
 \label{rem showed and existings(a)}
 When $m=0$, the bounds $S_{0,n}(g)$ are the bounds $S_n$
 given by \citet[Theorem 4.1, pp.~179--180]{APP1}, see
 \eqref{eq existing Poincare};
 \item
 \label{rem showed and existings(b)}
 The results of \Cref{theo variance bounds}
 also apply to the special case when $n=0$ (note that the second
 sum is empty and is treated as zero). In this case, the
 lower bound $S_{m,0}(g)$ is reduced to the one given by
 \citet[Theorem 4.1, pp.~518--519]{APP2}, see \eqref{eq
 Variance bounds Bessel-type}.
 \end{enumerate}
 \end{remark}

 \begin{remark}
 \label{rem infty bounds}
 Regarding the conditions of \Cref{theo variance bounds} imposed on the function $g$,
 we note that
 $g\in\CH^{\max\{m,n\},n-1}(X)\smallsetminus\CH^{\max\{m,n\},n}(X)$
 implies that the bound $S_{m,n}(g)$ is trivial, i.e.,
 $+\infty$ when $n$ is odd and $-\infty$ when $n$ is even.
 Of course, such a $g$ exists only when the
 support is infinite (with $\delta=0$).
 \end{remark}

 When $M<\infty$ and $m+n=M$, then $R_{m,n}(g)=0$ and the
 variance bound $S_{m,n}(g)$ is equal to $\Var[g(X)]$
 for any $g$. In any other case, it is of some interest to find
 an upper bound for the residual $R_{m,n}(g)$.

 \begin{proposition}\label{prop bounds of residual}
 Assume the conditions of \Cref{theo variance bounds}, with $m+n<M$, and, further, suppose that
 $g\in\CH^{T,T}(X)$ for some $T\in\{n,\ldots,m+n+1\}$.
 Then,
 the residual $R_{m,n}(g)$, given by \eqref{eq R_m:n},
 is bounded above by
 \begin{equation}
 \label{eq E(g^(tau)(X))^2}
 u_{\tau}\E\left\{{q^{[\tau]}(X)}\left[\Delta^{\tau}g(X)\right]^2\right\}, \quad \tau=n,n+1,\ldots,T,
 \end{equation}
 where
 $u_{\tau}
 =u_{m,n,\tau}(X)\doteq{\Pi_\delta^{[n]}(2m+n+1)}\left/\left\{{{m+n\choose{n}}(m+n+1)_\tau\Pi_\delta^{[n+\tau]}(m)}\right\}\right.$.
 \end{proposition}
 \begin{proof}
 Using \eqref{eq E q^[i](Delta^i g)^2}, we write the quantity
 \eqref{eq E(g^(tau)(X))^2} in the form
 $\sum_{k=\tau}^{M}\pi_{k;\tau}\alpha_k^2$. Next, consider the
 sequence
 $\{w_{k;\tau}={\pi_{k;\tau}}/{r_{k;m,n}(\delta)}\}_{k=m+n+1}^{M}$,
 where the numbers $r_{k;m,n}(\delta)$ are given by \eqref{eq
 R_m:n}, and observe that this sequence is increasing in
 $k$, with $w_{m+n+1;\tau}=1$.
 \end{proof}

 In general, the upper bounds (when there are at least two)
 of the residual $R_{m,n}(g)$, given by \eqref{eq
 E(g^(tau)(X))^2}, are not comparable.

 Next, for $n$ fixed, we investigate the bounds $S_{m,n}(g)$
 as $m$ increases.

 \begin{theorem}\label{theo S_n:m1 vs S_n:m2}
 Suppose $X\in\CX$ and fix a positive integer $n$ and a
 function $g\in\CH^{\tilde{m},n}(X)$, where
 $\tilde{m}$ (with $\tilde{m}\ge{n}$) can be finite
 or infinite. Then, for each $m_1,m_2$ such that $0 \le m_1 <
 m_2 \le \min\{\tilde{m},M\}$, the following inequality
 holds:
 \begin{equation}
 \label{eq R_n:m1 vs R_n:m2}
 |\Var[g(X)]-{S}_{m_1,n}(g)|\ge\zeta_{m_1,m_2,n}(\mu;q)|\Var[g(X)]-{S}_{m_2,n}(g)|,
 \end{equation}
 where $\zeta_{m_1,m_2,n}(\mu;q)=\zeta_{m_1,m_2,n}$
 is given by
 \begin{equation}
 \label{eq zeta_m1:m2:n}
 \zeta_{m_1,m_2,n}\doteq\left\{\begin{array}{@{}c@{\quad\textrm{if }}l@{}}
                    \frac{(m_2+n)_n(M-m_1-1)_n\Pi_\delta^{[n]}(m_2)\Pi_\delta^{[n]}(m_1+M)}
                         {(m_1+n)_n(M-m_2-1)_n\Pi_\delta^{[n]}(m_1)\Pi_\delta^{[n]}(m_2+M)},
                    & |M|<\infty,\\
                    \frac{(m_2+n)_n\Pi_\delta^{[n]}(m_2)}{(m_1+n)_n\Pi_\delta^{[n]}(m_1)},
                    & |M|=\infty.
                    \end{array}\right.
 \end{equation}
 For both cases, $|M|<\infty$ and $|M|=\infty$,
 \begin{equation}
 \label{eq bound for zeta}
 \zeta_{m_1,m_2,n}>{(m_2+n)_n}/{(m_1+n)_n}.
 \end{equation}
 The equality in \eqref{eq R_n:m1 vs R_n:m2} holds if
 and only if the function $g\colon S\to\RR$ is identically equal
 to a
 polynomial of degree at most $n+m_1$.
 \end{theorem}
 \begin{proof}
 Note that if $n+m_2=M$, then $S_{m_2,n}(g)=\Var[g(X)]$ for
 every function $g$ and \eqref{eq R_n:m1 vs R_n:m2} holds
 in a trivial way. Otherwise, we consider the finite or
 infinite positive sequence
 \[
 \left\{\zeta_{k}=\frac{r_{k;m_1,n}(\delta)}{r_{k;m_2,n}(\delta)}
                  =\frac{(m_2+n)_n(k-m_1-1)_n\Pi_\delta^{[n]}(m_2)\Pi_\delta^{[n]}(m_1+k)}
                        {(m_1+n)_n(k-m_2-1)_n\Pi_\delta^{[n]}(m_1)\Pi_\delta^{[n]}(m_2+k)}
  \right\}_{k=m_2+n+1}^{M}.
  \]

 \begin{quote}
 {\normalsize
 {\bfseries Claim.}
  \it
 The sequence $\{\zeta_{k}\}_{k=m_2+n+1}^{M}$ is strictly
 decreasing in $k$.
 \smallskip

 \noindent
 {\bfseries Proof of Claim.}
 \rm
 Since $\{{r_{k;m_1,n}(\delta)}/{r_{k+1;m_1,n}(\delta)}\}/\{{r_{k;m_2,n}(\delta)}/{r_{k+1;m_2,n}(\delta)}\}={\zeta_{k}}/{\zeta_{k+1}}$,
 $k=m_2+n+1,\ldots,M-1$, it is sufficient to show that the function
 $h(m)={r_{k;m,n}(\delta)}/{r_{k+1;m,n}(\delta)}
 ={(k-m-n)[1-(m+k)\delta]}/\{(k-m)[1-(m+n+k)\delta]\}$,
 $0\le m \le M-n-1$,
 is strictly decreasing. After some algebra,
 $h'(m)=-{n[1-(2m+n)\delta](1-2k\delta)}\left/\left\{(k-m)^2[1-(m+n+k)\delta]^2\right\}\right.$.
 If $\delta\le0$, then it is obvious that $h'(m)<0$; if
 $\delta>0$, then it is necessary that $M<\infty$ and, using
 \Cref{rem finite S: delta>0}, again it follows that
 $h'(m)<0$ and the
 claim is  proved.\hfill$\square$}
 \end{quote}

 If $M<\infty$, then the Claim shows that
 $\min_{k\in\{m_2+n+1,\ldots,M\}}\{\zeta_k\}=\zeta_M=\zeta_{m_1,m_2,n}$.
 If $M=\infty$, then observe that
 \begin{equation}
 \label{eq lim zeta_k}
 \zeta_k\searrow\frac{(m_2+n)_n\Pi_\delta^{[n]}(m_2)}{(m_1+n)_n\Pi_\delta^{[n]}(m_1)}=\zeta_{m_1,m_2,n}
 \quad\textrm{as } k\to\infty.
 \end{equation}
 Moreover, observing that ${r_{k;m_1,n}(\delta)}>0$ and
 ${r_{k;m_2,n}(\delta)}=0$ for all $k=n+m_1+1,\ldots,n+m_2$,
 \eqref{eq R_n:m1 vs R_n:m2} follows.

 If $\delta=0$ and $M=\infty$, then \eqref{eq bound for
 zeta} is obvious. For $\delta\le0$ and $M<\infty$, we observe that
 $\zeta_M>{(m_2+n)_n\Pi_\delta^{[n]}(m_2)}\left/\left[(m_1+n)_n\Pi_\delta^{[n]}(m_1)\right]\right.$,
 see \eqref{eq lim zeta_k},
 and \eqref{eq bound for zeta} follows.
 Now, assume $\delta>0$
 ($M<\infty$). Since
 ${\Pi_\delta^{[n+M-k]}(m_1+k)}>{\Pi_\delta^{[n+M-k]}(m_2+k)}>0$,
 it is sufficient to show that
 ${(M-m_1-1)_n\Pi_\delta^{[n]}(m_2)}\ge{(M-m_2-1)_n\Pi_\delta^{[n]}(m_1)}>0$.
 Observing that
 ${(M-m_1-1)_n\Pi_\delta^{[n]}(m_2)}\left/\left\{(M-m_2-1)_n\Pi_\delta^{[n]}(m_1)\right\}\right.
 =\prod_{j=0}^{n-1}\frac{(M-n+j-m_1)_n[1-(m_2+j)\delta]}{(M-n+j-m_2)_n[1-(m_1+j)\delta]}$,
 and putting $\eta_j\mapsto M-n+j$ and $\xi_j\mapsto
 1-j\delta$,
 it is sufficient to show that
 $[(\eta_j-m_1)(\xi_j-m_2)]/[(\eta_j-m_2)(\xi_j-m_1)]>1$
 for all $j=0,\ldots,n-1$.
 This is equivalent to
 $\xi_j-\eta_j\delta>0$, that is, $\delta<(M-n+2j)^{-1}$
 for all $j=0,\ldots,n-1$. Observe that
 for each $j=0,\ldots,n-1$,
 $(M-n+2j)^{-1}\ge[M-n+2(n-1)]^{-1}=(M+n-2)^{-1}\ge(2M-2)^{-1}=[2(|S|-2)]^{-1}>\delta$;
 see \Cref{rem finite S: delta>0}. Thus, \eqref{eq zeta_m1:m2:n} holds in any case. Finally, writing
 $|\Var[g(X)]-{S}_{m_1,n}(g)|-\zeta_{m_1,m_2,n}
 |\Var[g(X)]-{S}_{m_2,n}(g)|=\sum_{k=n+m_1+1}^{M}\theta_{k}
 \alpha_k^2$,
 we observe that $\theta_{k}>0$ for all $k$. Thus, the
 equality in \eqref{eq R_n:m1 vs R_n:m2} holds if and only
 if $g$ is identified with a polynomial of degree at most
 $n+m_1$.
 \end{proof}

 \begin{remark}
 \label{rem infty2}
 Assume the conditions of \Cref{theo S_n:m1 vs S_n:m2}.
 \begin{enumerate}[itemsep=.5ex, wide, labelwidth=!, labelindent=0pt, label=\rm(\alph*), ref=\textcolor{black}{\alph*}]
 \item
 \label{rem infty2(a)}
 In view of \Cref{rem S_m:n with n fixed}\eqref{rem S_m:n with n fixed(a)}, the
 bounds $\{S_{m,n}(g)\}_{m=0}^{\tilde{m}}$ are of the
 same kind. From \eqref{eq R_n:m1 vs R_n:m2}, it follows
 that the bound ${S}_{m_2,n}(g)$ is better than the bound
 ${S}_{m_1,n}(g)$. Thus, writing $n=2r$ (when $n$ is even)
 and $n=2r+1$ (when $n$ is odd), we have
 \[
 S_{0,2r}(g)\le S_{1,2r}(g)\le\cdots\le\Var[g(X)]\le\cdots\le S_{1,2r+1}(g)\le S_{0,2r+1}(g);
 \]
 \item
 \label{rem infty2(b)}
 For the case $\tilde{m}=M=\infty$, from \eqref{eq
 Var g(X)=sum_k=0^M a_k^2}, \eqref{eq S_m:n} and \eqref{rem infty2(a)}, it
 follows that
  \[
  \substack{\displaystyle S_{m,n}(g)\nearrow\Var[g(X)] \\ \mbox{\footnotesize[when $n$ is even]}}
  \quad \substack{\displaystyle \textrm{or} \\ \ } \quad
  \substack{\displaystyle S_{m,n}(g)\searrow\Var[g(X)],\\ \mbox{\footnotesize[when $n$ is odd]}}
  \quad
  \substack{\displaystyle \textrm{as}\quad m\to\infty. \\ \ }
  \]
 \end{enumerate}
 \end{remark}

 Now, we compare the existing variance bound $S_{0,n}(g)$,
 see \Cref{rem showed and existings}\eqref{rem showed and existings(a)}, with the best
 proposed bound shown in this section, requiring the same
 conditions on $g$, i.e., with the bound $S_{n,n}(g)$, see
 \Cref{rem S_m:n with n fixed}\eqref{rem S_m:n with n fixed(b)}.
 \begin{corollary}\label{cor S_{n:n} vs S_{0:n}}
 The variance bounds $S_{n,n}(g)$ and $S_{0,n}(g)$ are of
 the same kind and require the same assumptions on $g$.
 Moreover, the new bound $S_{n,n}(g)$ is better than the
 existing (see \Cref{rem showed and existings})
 bound $S_{0,n}(g)$. Specifically,
 \[
 |\Var[g(X)]-S_{0,n}(g)|\ge\zeta_{0,n,n}|\Var[g(X)]-{S}_{n,n}(g)|,
 \]
 with $\zeta_{0,n,n}>{2n\choose{n}}$. The equality holds
 only in the trivial case when
 $\Var[g(X)]=S_{n,n}(g)=S_{0,n}(g)$, i.e., the function
 $g\colon S\to\RR$ is identified with a polynomial of degree at most
 $n$.
 \end{corollary}

  \begin{remark}
 \label{remUVMUe}
 Assume that $X_1,\ldots,X_\nu$ is a random sample from the geometric
 distribution with parameter $\theta\in(0,1)$,
 i.e., with pmf $p(j)=\theta(1-\theta)^j$, $j=0,1,\ldots$,
 and let $X=X_1+\cdots+X_\nu$ be the complete sufficient statistic.
 The uniformly minimum variance unbiased estimator of $-\log(\theta)$
 is $T_\nu=T_\nu(X)=\sum_{j=\nu}^{\nu+X-1}1/j$.
 Variance bounds of the kind of \Cref{theo variance bounds}
 have been used for constructing bounds of $\Var(T_\nu)$;
 see \citet[Section~5]{APP1} and \citet[Application~5.1]{APP2}.
 In the similar and easy manner, we can use the results of
 \Cref{theo variance bounds,theo S_n:m1 vs S_n:m2}
 in regard to the approximation of $\Var(T_\nu)$
 and its accuracy.
 \end{remark}

 \noindent
 {\bf Acknowledgements.} This research
 has been co-financed by the European Union (European Social
 Fund -- ESF) and Greek national funds through the Operational
 Program ``Education and Lifelong Learning'' of the National
 Strategic Reference Framework (NSRF) -- Research Funding
 Program: ARISTEIA, Grant No.: 4357. Also, this
 work is partially
 supported by the University of Athens Research Grant
 70/4/5637 and by internal
 funds, Department of Biostatistics, SUNY Buffalo.
 This work was also partially supported by the Natural
 Sciences and Engineering Research Council of Canada through
 an Individual Discovery Grant to the second author.
 Furthermore, the authors acknowledge the editorial team
 who handled the paper for providing suggestions that
 resulted in improving the presentation of the results.

\end{document}